\documentclass[preprint,authoryear]{elsarticle}

\biboptions{sort&compress}

\usepackage{graphicx}      % include this line if your document contains figures
\usepackage{natbib}        % required for bibliography
\usepackage{amssymb}
\usepackage{amsmath}
\usepackage{enumerate}
\usepackage{placeins}
\usepackage{epstopdf}
\usepackage[none]{hyphenat}
\newtheorem{thm}{Theorem}

\newtheorem{corol}{Corollary}

\newproof{pf}{Proof}

\begin{document}
\bibliographystyle{model5-names}
\pagestyle{plain}
\setcounter{page}{0}
\pagenumbering{arabic}
\begin{frontmatter}

\title{Lipschitz Constants in Experimental Optimization} 

\author[ind]{Gene A. Bunin\corref{cor1}}
\ead{gene.bunin@ronininstitute.org}
\author[edin]{Gr\'egory Fran\c cois}
\ead{gregory.francois@ed.ac.uk}

\cortext[cor1] {Author to whom correspondence should be addressed (tel.: +86 13201322405).}

\address[ind]{Ronin Institute for Independent Scholarship}
\address[edin]{Institute for Materials and Processes, School of Engineering, The University of Edinburgh, Edinburgh EH9 3FB, UK}

\begin{abstract}                % Abstract of not more than 250 words.

The Lipschitz constant of a response surface function upper bounds the sensitivity of a dependent variable to changes in the independent ones. Traditionally, such constants have found much implicit and abstract use in mathematically oriented applications, but their potential for explicit use in more engineering-based domains has not been explored. The latter point is the subject of this paper, where we propose several ways in which the Lipschitz constants may be used explicitly in the domain of \emph{experimental optimization}. Specifically, we focus on how they may help ensure the satisfaction of constraints and on their potential role in reducing the negative effects of measurement or estimation uncertainty. A number of refinements to the proposed approaches are also derived, and some techniques for setting the constants are presented. 

\;


\end{abstract}

\end{frontmatter}
%===============================================================================

\section{Introduction: Experimental Optimization and Lipschitz Constants}
\label{sec:intro}

Let us consider the function $f_p : \mathbb{R}^{n_u} \rightarrow \mathbb{R}$ to be an \emph{experimental} relationship between the independent-variable vector, ${\bf u} \in \mathbb{R}^{n_u}$, and an experimental output quantity, $y$:

$$
y = f_p({\bf u}).
$$

\noindent Here, the word ``experimental'' will imply that evaluating the corresponding $y$ for a given ${\bf u}$ requires carrying out a \emph{physical} experiment that \emph{cannot be performed solely with a computer}, with the evaluation intrinsic to $f_p$ defined as a physical act and not as a series of numerical computations.

Such relationships appear virtually everywhere in the sciences \citep{Montgomery2012}, and are typically accompanied by systematic studies that attempt to discover the nature of the given $f_p$. Among the documented qualitative work on the subject, it is probably John Lind's 1747 examination of scurvy patients \citep{Dunn1997} that stands out as the earliest. On the quantitative front, it is likely to Ronald Fisher and his work on the design of experiments \citep{Fisher1935} that modern scientific investigation is most indebted, with the foundations for obtaining a mathematical approximation of $f_p$ laid out therein still very much in use today.

In many cases, however, simply identifying an experimental relationship is not the end goal. Instead, one is often searching to manipulate experimental conditions in such a way so as to obtain the ``best'' response while satisfying a number of safety or physical restrictions. In such scenarios, the \emph{experimental optimization} problem

\begin{equation}\label{eq:mainprob}
\begin{array}{rll}
\mathop {{\rm{minimize}}}\limits_{\bf{u}} & \phi_p ({\bf{u}}) &  \\
{\rm{subject}}\hspace{1mm}{\rm{to}} & g_{p,j}({\bf{u}}) \leq 0, & j = 1,...,n_{g_p} \vspace{1mm} \\
 & g_{j}({\bf{u}}) \leq 0, & j = 1,...,n_{g} \vspace{1mm} \\
& u_i^L \leq u_i \leq u_i^U, & i = 1,...,n_u
\end{array}
\end{equation}

\noindent usually arises, where $\phi : \mathbb{R}^{n_u} \rightarrow \mathbb{R}$ denote the cost function to be minimized, the functions $g : \mathbb{R}^{n_u} \rightarrow \mathbb{R}$ denote the $n_{g_p} + n_g$ constraints, and the constants $u_i^L, u_i^U$ denote the lower and upper limits on ${\bf u}$ (now termed ``decision variables''), respectively. The subscript $p$ is used to denote explicitly those relationships that are experimental in nature, its absence indicating that the function is numerical and only requires a computer or basic algebra to evaluate. A diverse range of practical engineering problems may be cast and solved in the form of (\ref{eq:mainprob}) \citep{ExpOpt}.

Because the experimental functions $\phi_p$ and $g_{p,j}$ are unknown and not subject to numerical evaluation, solving (\ref{eq:mainprob}) entails running a \emph{series} of experiments at ${\bf u}_0, {\bf u}_1, \hdots, {\bf u}_k$, where each new experiment is chosen via some algorithmic function, $\Gamma$, of the previously applied ${\bf u}$ and the corresponding experimentally evaluated function values:

\begin{equation}\label{eq:bigalgo}
{\bf u}_{k+1} = \Gamma \left( \begin{array}{l} {\bf u}_0, {\bf u}_1, \hdots, {\bf u}_k, \phi_p ({\bf u}_0), \phi_p ({\bf u}_1), \hdots, \phi_p ({\bf u}_k), \\ g_{p,j} ({\bf u}_0), g_{p,j} ({\bf u}_1), \hdots, g_{p,j} ({\bf u}_k), \; j = 1,...,n_{g_p} \end{array}  \right),
\end{equation}

\noindent the index $k$ used to denote the latest experimental iteration. The seminal paper with regard to solving (\ref{eq:mainprob}) is due to \cite{Box:51}, where $\Gamma$ essentially takes the form of a standard gradient-based descent method. Numerous other methods ($\Gamma$) and frameworks have been proposed in various mathematical and engineering domains since \cite[\S 2]{Bunin2013SIAM}.

In the present work, we explore a single facet of Problem (\ref{eq:mainprob}) and focus on the case where the experimental functions involved are \emph{Lipschitz continuous} over the explored domain, as we believe that making this assumption may benefit the solution procedure in multiple ways. Given the Lipschitz continuity of a function, we first state the fundamental law \citep[\S 2.2]{Biegler2010} that there exists a finite \emph{Lipschitz constant}, $\kappa$, such that

\begin{equation}\label{eq:lipgen}
f_p({\bf u}_b) \leq f_p ({\bf u}_a) + \kappa \| {\bf u}_b - {\bf u}_a \|_2, \;\; \forall {\bf u}_a, {\bf u}_b \in \mathcal{I},
\end{equation}

\noindent where $\mathcal{I} = \{ {\bf u} : u_i^L \leq u_i \leq u_i^U, \; i = 1,...,n_u \}$ is the \emph{experimental space} to which the experiments are restricted. As with any assumption, the practical validity of (\ref{eq:lipgen}) must be questioned, with the authors' empirical experience so far suggesting that the assumption is not unreasonable for many, if not most, practical problems. In these problems, the constant $\kappa$ may be understood as a maximal \emph{sensitivity} of the experimental function $f_p$ to the changes in the decision variables ${\bf u}$ over the experimental space $\mathcal{I}$.

Let us suppose now that an explicit value of $\kappa$ satisfying (\ref{eq:lipgen}) is \emph{known} or, at the very least, that a reliable overestimate is available. As will be shown in Section \ref{sec:uses}, such information opens several new doors, thereby making it possible to (a) satisfy important experimental constraints during the optimization process, and (b) reduce the potentially harmful effects of measurement or estimation uncertainty. Of the two, (b) is seen as useful while (a) is considered crucial, as it allows the introduction of \emph{simple} safety guarantees for problems that traditionally lack them. A number of simulated case-study examples will be presented to illustrate the potential benefits of the techniques.

Following the presentation of these core ideas, we will explore, in Section \ref{sec:refine}, the use of additional assumptions to derive an alternate version of (\ref{eq:lipgen}) that is, in most cases, tighter and thus less conservative. Specifically, it will be shown that 

\begin{itemize}
\item using multiple Lipschitz constants for $f_p$ (sensitivities with respect to individual decision variables), 
\item using local values rather than those valid for \emph{all} ${\bf u}_a$ and ${\bf u}_b$ in $\mathcal{I}$, and 
\item taking into consideration the potential convexity and concavity properties of $f_p$ 
\end{itemize}

\noindent all allow for a version of (\ref{eq:lipgen}) that, though visually more complex, can often lead to significant gains in performance. As may be expected, the price to pay for the improvement is that of additional assumptions. However, they are not compulsory, and the user may employ as many as can reasonably be made. 

Finally, Section \ref{sec:estim} of the paper will focus on the fundamental and important issue of setting the $\kappa$ values, which is not so trivial to do well -- here, one must remember that since $f_p$ is unknown, so is $\kappa$. A number of various methods that combine physical principles, model-based estimation, and data-driven refinement are discussed, with a summary of the authors' own empirical experience with these methods' effectiveness concluding the section.

\section{The Different Uses of the Lipschitz Constants}
\label{sec:uses}

Exploiting the Lipschitz constants will, in general, be synonymous with exploiting the corresponding lower and upper \emph{Lipschitz bounds},

\begin{equation}\label{eq:lipgenLU}
f_p ({\bf u}_a) - \kappa \| {\bf u}_b - {\bf u}_a \|_2 \leq f_p({\bf u}_b) \leq f_p ({\bf u}_a) + \kappa \| {\bf u}_b - {\bf u}_a \|_2, \;\; \forall {\bf u}_a, {\bf u}_b \in \mathcal{I},
\end{equation}

\noindent of which the lower bound is a symmetrical result easily obtained by switching ${\bf u}_a \leftrightarrow {\bf u}_b$ in (\ref{eq:lipgen}). Given the experimental nature of $f_p$, the practical implications of (\ref{eq:lipgenLU}) are important as they allow us to analyze the worst-case behavior of the function with respect to its departure from ${\bf u}_a$. In other words, if the experiment at ${\bf u}_a$ has been conducted but the one at ${\bf u}_b$ has not, it becomes possible to bound the greatest possible change in $f_p$ if ${\bf u}_b$ were applied.

This section focuses on the three major uses of (\ref{eq:lipgenLU}), outlining their theoretical foundations, discussing their implementation, and providing illustrations via several case-study examples.  

\subsection{Satisfaction of Experimental Constraints During Optimization}
\label{sec:consat}

In general, there exists no guarantee that the iterates ${\bf u}_1, {\bf u}_2,\hdots,{\bf u}_{k+1}$ generated by the algorithm $\Gamma$ of (\ref{eq:bigalgo}) will satisfy the experimental constraints $g_{p,j} ({\bf u}_{\bar k}) \leq 0$ for all $\bar k = 1,...,k+1$. The reason is simple -- the functions are unknown to the user, and this makes it impossible to know in advance if the next applied set of decision variables, ${\bf u}_{k+1}$, will satisfy $g_{p,j} ({\bf u}_{k+1}) \leq 0$. While different methodologies for solving (\ref{eq:mainprob}) do employ different techniques for ensuring that these constraints are met in some sense, the current state of the art is not very developed in this regard, with rigorous constraint satisfaction a notable gap in the experimental-optimization literature. When looking at the three prominent domains concerned with solving Problem (\ref{eq:mainprob}), we notice the following:

\begin{itemize}
\item derivative-free optimization methods, developed mostly by the mathematical community and often dealing with numerical optimization problems, only take precautions against constraints \emph{asymptotically} -- i.e., they ensure, usually via penalty-function methods, that ${\bf u}_{k+1}$ satisfies the constraints as $k \rightarrow \infty$ \citep{Liuzzi2010,Biegler2014} \citep[\S 14]{Conn2000};
\item response-surface optimization methods often deal with experimental systems and are used in a wide number of fields \citep{Myers2009}, but neither do they address the problem of constraint satisfaction during the experiments used to construct the response-surface model, nor do they rigorously attack the potential issue of the model optimum violating the constraints \citep{Umland1959,Michaels1963};
\item real-time optimization, as typically done in the chemical-engineering community \citep{Chachuat2009,Darby2011}, is arguably the most concerned with consistent constraint satisfaction, since $g_{p,j} ({\bf u}) > 0$ may represent a dangerous or economically disastrous operating regime; however, while both theoretical \citep{Loeblein1999,Zhang2002} and practical \citep{Govatsmack2005} methods have been developed for promoting the satisfaction of constraints, the former are often too complex and restrictive while the latter suffer from a lack of rigorous guarantees, with both types of methods tending to find solutions that are significantly suboptimal because of conservatism \citep{Li2008,Quelhas:12}.
\end{itemize}

\noindent In particular, for those problems where rigorous experimental-constraint satisfaction is highly desired, the above methods are largely inadequate, and only in the real-time optimization community is the rigorous satisfaction of constraints at \emph{every} experiment even a concern.

Some years ago, the authors proposed a Lipschitz approach to this problem \citep{Bunin2011}. Striking what we believe to be a good balance between rigor and simplicity, the method suggests to start with an initially feasible ${\bf u}_0$, where $g_{p,j} ({\bf u}_0) \leq 0$, and to repeatedly couple this with the upper Lipschitz bound of (\ref{eq:lipgenLU}) to generate a series of conditions that guarantee constraint satisfaction throughout:

\begin{equation}\label{eq:lipseq}
\begin{array}{l}
\displaystyle g_{p,j} ({\bf u}_{0}) + \kappa_{g_p,j} \| {\bf u}_{1} - {\bf u}_{0} \|_2 \leq 0 \Rightarrow g_{p,j} ({\bf u}_1) \leq 0 \vspace{1mm} \\
\displaystyle g_{p,j} ({\bf u}_{1}) + \kappa_{g_p,j} \| {\bf u}_{2} - {\bf u}_{1} \|_2 \leq 0 \Rightarrow g_{p,j} ({\bf u}_2) \leq 0 \vspace{1mm} \\
\vdots \vspace{1mm} \\
\displaystyle g_{p,j} ({\bf u}_{k}) + \kappa_{g_p,j} \| {\bf u}_{k+1} - {\bf u}_{k} \|_2 \leq 0 \Rightarrow g_{p,j} ({\bf u}_{k+1}) \leq 0,
\end{array}
\end{equation}

\noindent where $\kappa_{g_p,j}$ denotes the corresponding Lipschitz constant for $g_{p,j}$.

The algorithmic implementation of (\ref{eq:lipseq}) is straightforward, as at the current iterate $k$ the value of $g_{p,j}({\bf u}_k)$ has been obtained, $\kappa_{g_p,j}$ is set, ${\bf u}_k$ has been applied, and ${\bf u}_{k+1}$ is controlled by the user. As such, it suffices to add the constraint

\begin{equation}\label{eq:congamma}
g_{p,j}({\bf u}_k) + \kappa_{g_p,j} \| {\bf u}_{k+1} - {\bf u}_k \|_2 \leq 0
\end{equation}

\noindent to $\Gamma$ when selecting the next ${\bf u}_{k+1}$ to apply. Clearly, because $g_{p,j}({\bf u}_k) \leq 0$, one can always bring ${\bf u}_{k+1}$ sufficiently close to ${\bf u}_k$ so as to ensure that (\ref{eq:congamma}) is satisfied. The simplicity of the implementation is particularly attractive since no complex model of the experimental functions is required, which is the case in some of the rigorous constraint-satisfaction methods previously proposed \citep{Loeblein1999,Zhang2002}.

The rigor of (\ref{eq:lipseq}) should also be apparent. If $g_{p,j} ({\bf u}_0) \leq 0$ and $\kappa_{g_p,j}$ is such that (\ref{eq:lipgen}) is satisfied, then (\ref{eq:lipseq}) guarantees that the sequence of ${\bf u}_1, {\bf u}_2,\hdots,{\bf u}_{k+1}$  meets the constraints at every experiment. Practically speaking, this is a very useful result.

Of course, this approach is not without its caveats, and it would be dangerous to assume otherwise. The most major is the correctness of $\kappa_{g_p,j}$, as setting $\kappa_{g_p,j}$ too small essentially invalidates the method's rigor and makes its guarantees unreliable. The obvious solution of making $\kappa_{g_p,j}$ ``very big'' comes with the major drawback of conservatism, and so cannot be recommended as it would implicitly shrink the distances between ${\bf u}_2$ and ${\bf u}_1$, ${\bf u}_3$ and ${\bf u}_2$, and so on, resulting in the optimization advancing very slowly. Given the importance of this issue, we will return to it in detail in Section \ref{sec:estim}, assuming until then that a correct choice of $\kappa_{g_p,j}$ is available.

\begin{figure}
\begin{center}
\includegraphics[width=6cm]{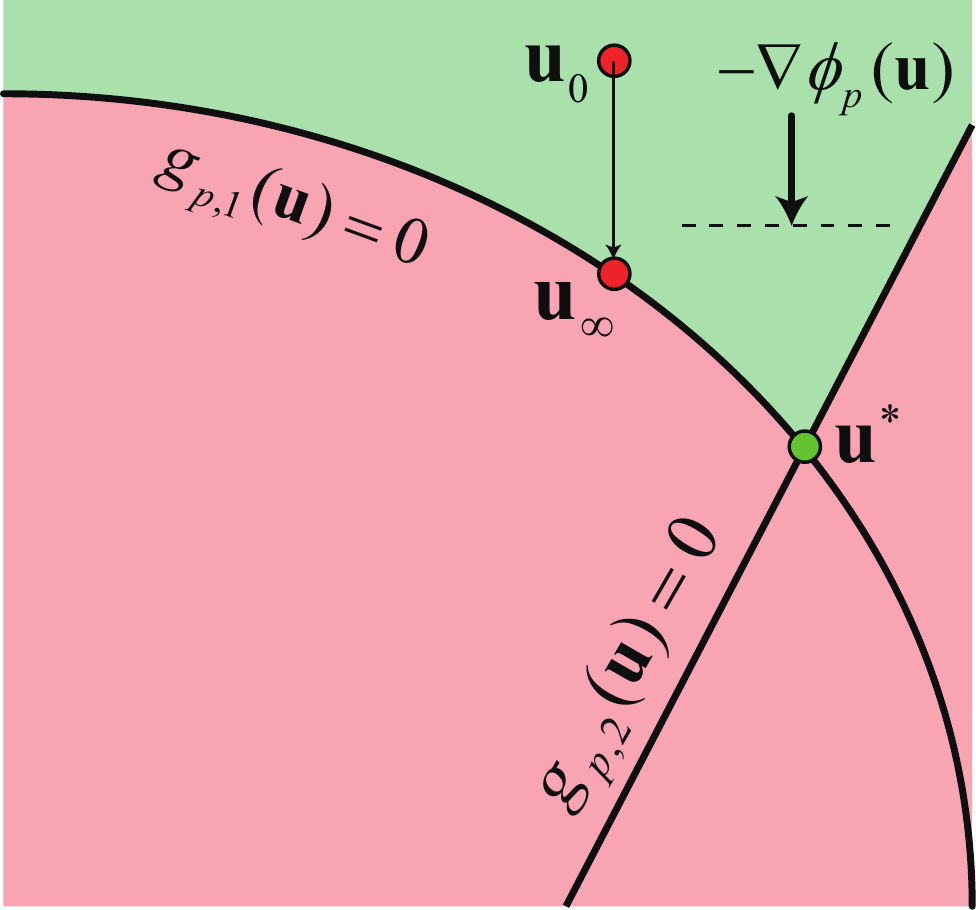}
\caption{An example of premature convergence for a problem with a linear cost function where the applied $\Gamma$ continuously takes steps in the gradient descent direction while keeping the constraint satisfied via the Lipschitz condition (\ref{eq:congamma}). Gradually, the algorithm converges to the point ${\bf u}_\infty$ on the boundary of an experimental constraint ($g_{p,1}$), after which it can no longer make any progress towards the true optimum ${\bf u}^*$, as any step chosen would violate (\ref{eq:congamma}). The green and red areas represent the feasible (all constraints satisfied) and infeasible (some constraints violated) sets of the problem, respectively.}
\label{fig:feasbreakSIAM}
\end{center}
\end{figure}

A lesser but important concern lies in the potential issue of \emph{premature convergence} that may be experienced by any $\Gamma$ that adds (\ref{eq:congamma}) to its definition. More concretely, this refers to the case where the experiments generated by $\Gamma$ approach an experimental constraint faster than they approach an optimum, as this drives a given $g_{p,j}$ to 0 and forces ${\bf u}_{k+1} \approx {\bf u}_{k}$ for (\ref{eq:congamma}) to be met, essentially forcing the algorithm to terminate on some active experimental constraint. We illustrate the issue geometrically in Fig. \ref{fig:feasbreakSIAM}. Although the solution to this particular problem will not be addressed in this paper, we do note that it may be attained with the addition of supplementary conditions that force $\Gamma$ to not approach any constraint ``too quickly''. The interested reader is referred to the unpublished document of \cite{Bunin2013SIAM} for the details.

Finally, the method does rely on the initial experiment satisfying the constraints. However, this is not expected to be a limitation, since it is typically assumed that the experimental optimization of some given system will start at a point that is safe but suboptimal. Indeed, if ${\bf u}_0$ were not safe, then one would not expect optimization to be the user's top priority.

\subsubsection{Case-Study Example: Constraint Adaptation to Minimize the Batch Time of Polystyrene Production}
\label{sec:poly}

To illustrate how the Lipschitz method may be successfully incorporated into the solution of an experimental optimization problem, let us take the case-study example of minimizing the batch time of a polystyrene production reactor, originally reported by \cite{Gentric1999}. Here, we will use the problem formulation of \cite{Francois2005}, where one is interested in choosing the ``switching times'' of the reactor's temperature profile so as to both (a) minimize the time required to reach the desired conversion, favored by higher temperatures, and (b) meet a terminal lower limit on the molecular weight of the product, favored by lower temperatures. The original problem is a dynamic optimization problem as it seeks to find an optimal \emph{profile}, but the piecewise definition of the profile via prescribed ``arcs'' and switching times essentially turns this into an experimental optimization problem, which can then be solved in a batch-to-batch manner by varying the switching times (and thus the profile) from batch to batch until an optimal profile is hopefully found. As formulated, the problem has an experimental cost function (the time of the batch) and a single experimental constraint (the lower limit on the molecular weight). Because the molecular weight is a product specification, it would be very wasteful to run batches that violate this constraint, as all such batches may be discarded as inadequate for failing to meet specifications. It follows that an effort to satisfy this constraint for all batches should be made.

The exact problem solved is Problem P3 from the ExpOpt database \citep{ExpOpt}, which the interested reader may access to find a more detailed problem description, the steps taken to place the problem into the standard form (\ref{eq:mainprob}), and the actual code used to simulate the case study. Only the problem as written in standard form,

\begin{equation}\label{eq:polyprob}
\begin{array}{rl}
\mathop {{\rm{minimize}}}\limits_{\bf{u}} & \phi_p ({\bf{u}}) := t_b ({\bf u}) \\
{\rm{subject}}\hspace{1mm}{\rm{to}} & g_{p,1}({\bf{u}}) := -M({\bf u}) + 2 \cdot 10^6 \leq 0 \vspace{1mm} \\
& 50 \leq u_1 \leq 450 \vspace{1mm} \\
& 600 \leq u_2 \leq 1000,
\end{array}
\end{equation}

\noindent is given here, where $t_b$ denotes the final batch time to be minimized (in seconds), $M$ denotes the number average molecular weight (in grams per mole), and $u_1$ and $u_2$ are the two switching times that define the temperature profile (in seconds). So as to avoid unnecessary numerical complications, (\ref{eq:polyprob}) is scaled prior to being solved, with $t_b$ divided by 8000, the constraint function divided by $4 \cdot 10^6$, and the two decision variables normalized to lie in the unit box defined by the coordinates $(0,0)$ and $(1,1)$. The initial experiment is chosen as ${\bf u}_0 := (242.39, 945.30)$.

To solve the problem, we apply as $\Gamma$ a constraint-adaptation \citep{Chachuat2008a}, or bias-update \citep{Forbes1994b}, algorithm. In applying this method, we suppose the availability of an approximate process model,

$$
\begin{array}{l}
\phi_{\hat p} ({\bf u}, \theta_1, \theta_2) \approx \phi_p ({\bf u}) \\
g_{\hat p, 1} ({\bf u},\theta_1, \theta_2) \approx g_{p,1} ({\bf u}),
\end{array}
$$

\noindent where $\theta_1$ and $\theta_2$ are the uncertain model parameters -- in this case, the rate constants for propagation and transfer to monomer, respectively. Because the values of the model, $\theta_1 = 4.6 \cdot 10^6$ and $\theta_2 = 2.3 \cdot 10^{11}$ (liters per mole per second), differ from the true values of $5.7 \cdot 10^6$ and $1.5 \cdot 10^{11}$, solving the model-based optimization problem

$$
\begin{array}{rl}
\mathop {{\rm{minimize}}}\limits_{\bf{u}} & \phi_{\hat p} ({\bf u},\theta_1,\theta_2) \\
{\rm{subject}}\hspace{1mm}{\rm{to}} & g_{\hat p,1}({\bf u},\theta_1,\theta_2) \leq 0 \vspace{1mm} \\
& 50 \leq u_1 \leq 450 \vspace{1mm} \\
& 600 \leq u_2 \leq 1000
\end{array}
$$

\noindent directly will not yield the optimal point. In constraint adaptation, one aims instead to iteratively approach a better solution by ensuring that the constraint values in the model-based optimization tend to the true measured values by applying a local bias correction term, $\epsilon_k$, and setting each next ${\bf u}_{k+1}$ as

\begin{equation}\label{eq:CAprob}
\begin{array}{rl}
{\bf u}_{k+1} := {\rm arg} \mathop {{\rm{minimize}}}\limits_{\bf{u}} & \phi_{\hat p} ({\bf u},\theta_1,\theta_2) \\
{\rm{subject}}\hspace{1mm}{\rm{to}} & g_{\hat p,1}({\bf u},\theta_1,\theta_2) + \epsilon_k \leq 0 \vspace{1mm} \\
& 50 \leq u_1 \leq 450 \vspace{1mm} \\
& 600 \leq u_2 \leq 1000,
\end{array}
\end{equation}

\noindent with $\epsilon_k$ defined as a weighted sum of the current model error and the previous correction term:

\begin{equation}\label{eq:CAupdate}
\epsilon_k := \alpha \left[ g_{p,1} ({\bf u}_k) - g_{\hat p,1} ({\bf u}_k, \theta_1, \theta_2) \right] + (1-\alpha) \epsilon_{k-1}. 
\end{equation}

\begin{figure}
\begin{center}
\includegraphics[width=8cm]{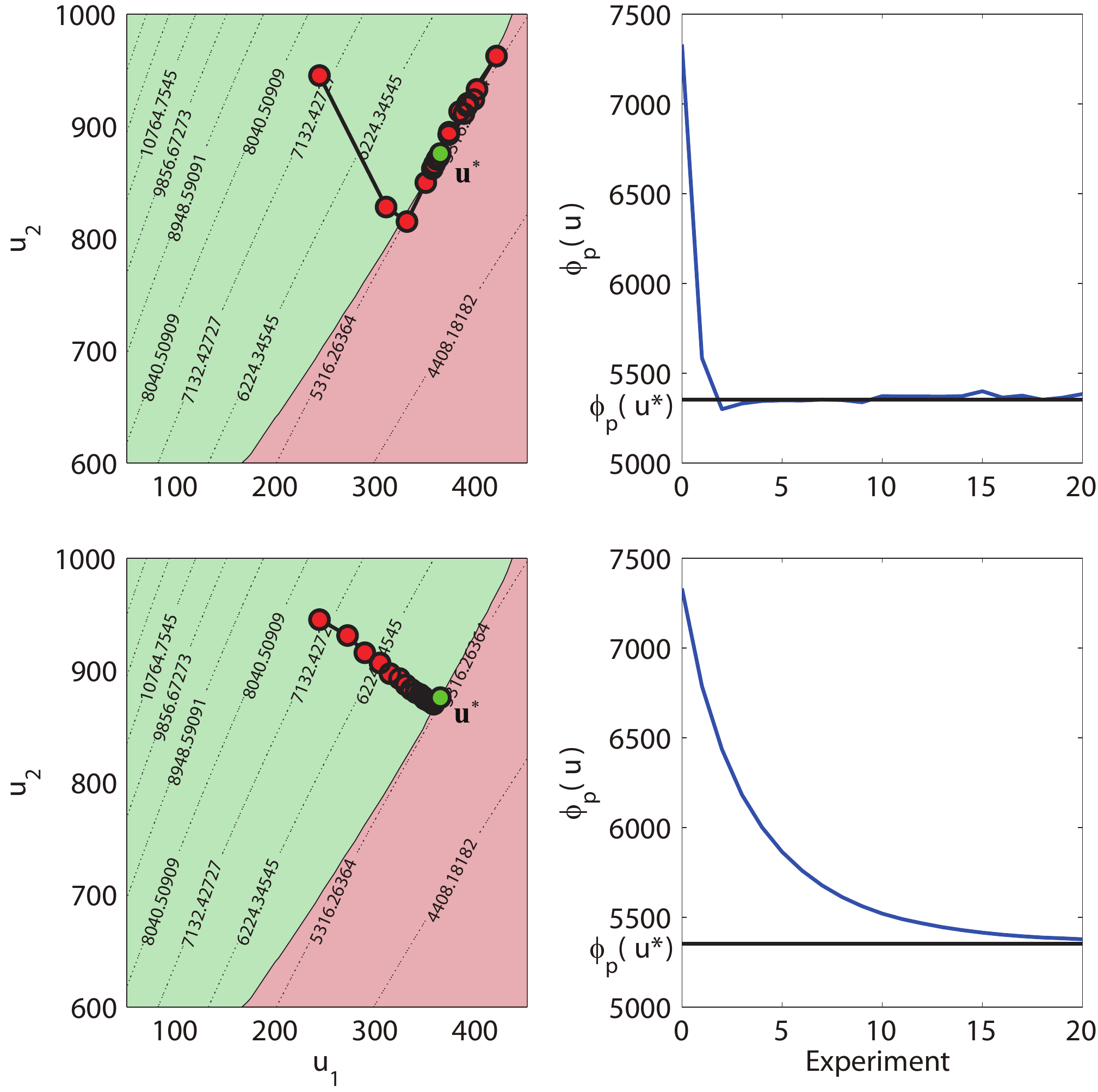}
\caption{The solution of Problem (\ref{eq:polyprob}) with the constraint-adaptation algorithms of (\ref{eq:CAprob}) (top) and (\ref{eq:CAproblip}) (bottom).}
\label{fig:CAex}
\end{center}
\end{figure}

\noindent In this example, we initialize $\epsilon_{-1}$ as 0, and apply a filter gain of $\alpha := 0.7$. Applying (\ref{eq:CAprob}) and (\ref{eq:CAupdate}) repeatedly essentially solves the problem, and it can be seen from the top set of results in Fig. \ref{fig:CAex} that the algorithm succeeds at finding the neighborhood of the solution fairly quickly. However, the path taken results in consistent violations of the molecular weight constraint.

To incorporate the Lipschitz bound, we simply add Condition (\ref{eq:congamma}) as a constraint in (\ref{eq:CAprob}):

\begin{equation}\label{eq:CAproblip}
\begin{array}{rl}
{\bf u}_{k+1} := {\rm arg} \mathop {{\rm{minimize}}}\limits_{\bf{u}} & \phi_{\hat p} ({\bf u},\theta_1,\theta_2) \\
{\rm{subject}}\hspace{1mm}{\rm{to}} & g_{\hat p,1}({\bf u},\theta_1,\theta_2) + \epsilon_k \leq 0 \vspace{1mm} \\
& 50 \leq u_1 \leq 450 \vspace{1mm} \\
& 600 \leq u_2 \leq 1000 \vspace{1mm}\\
& g_{p,1} ({\bf u}_k) + \kappa_{g_p,1} \| {\bf u} - {\bf u}_k \|_2 \leq 0,
\end{array}
\end{equation}

\noindent with a (correct) value\footnote{The actual issue of setting correct Lipschitz constants is addressed in more detail in Section \ref{sec:estim}. For the case-study examples in this section, we will, for simplicity, use correct Lipschitz constants, as may be obtained by exhaustive numerical tests on $f_p$, available to us in these simulated studies.} of 3 used for $\kappa_{g_p,1}$. The performance of the resulting $\Gamma$ is given as the bottom set of results in Fig. \ref{fig:CAex}, and one observes a trajectory of ${\bf u}$ that, though approaching the optimum slower, does so without \emph{any} violations.

As an aside, we note that this way of implementing (\ref{eq:congamma}) may be used for any $\Gamma$ that computes ${\bf u}_{k+1}$ by solving a numerical optimization problem. A more general alternative of employing (\ref{eq:congamma}) in a line search may be used otherwise \citep{Bunin2013SIAM}.  

\subsection{Satisfaction of Constraints During Local Perturbations}
\label{sec:excite}

The previous section focused on modifying the optimization algorithm $\Gamma$ so that the constraints remain satisfied at every future ${\bf u}_{k+1}$. However, $\Gamma$ may not be the only law that dictates what experiments are run, as it is often the case that \emph{additional experiments} are periodically required to gather information about the experimental function locally. For example, one may wish to perform experiments in the neighborhood of ${\bf u}_k$ so as to identify the parameters of a locally valid model \citep{Conn2000,Pfaff2006} or to estimate the local function derivatives \citep{Gao2005,Marchetti2010}, either of which may then be used to assist $\Gamma$ in computing ${\bf u}_{k+1}$. Naturally, there then arises a conflict between constraint satisfaction and obtaining information when one approaches a constraint, with the information-gathering perturbations becoming a potential constraint-violation hazard.

A viable solution to this conflict has already been discussed in detail by \cite{Bunin2016}, and what we present here is a simple variant. Namely, to guarantee that any point in the perturbation ball of $\mathcal{B}_e = \{ {\bf u} : \| {\bf u} - {\bf u}_k \|_2 \leq \delta_e \}$ satisfy a given constraint $g_{p,j} ({\bf u}) \leq 0$, it is sufficient to enforce that the value $g_{p,j} ({\bf u}_k)$ satisfy a constraint back-off of $\delta_e \kappa_{g_p,j}$:

\begin{equation}\label{eq:backoff}
g_{p,j} ({\bf u}_k) + \delta_e \kappa_{g_p,j}  \leq 0 \Rightarrow g_{p,j} ({\bf u}) \leq 0, \;\; \forall {\bf u} \in \mathcal{B}_e \cap \mathcal{I}.
\end{equation}

The validity of (\ref{eq:backoff}) is easily demonstrated. For every ${\bf u} \in \mathcal{B}_e$,

$$
\| {\bf u} - {\bf u}_k \|_2 \leq \delta_e \Rightarrow g_{p,j} ({\bf u}_k) + \kappa_{g_p,j} \| {\bf u} - {\bf u}_k \|_2  \leq g_{p,j} ({\bf u}_k) + \delta_e \kappa_{g_p,j},
$$

\noindent and thus

$$
g_{p,j} ({\bf u}_k) + \delta_e \kappa_{g_p,j} \leq 0 \Rightarrow g_{p,j} ({\bf u}_k) + \kappa_{g_p,j} \| {\bf u} - {\bf u}_k \|_2 \leq 0 \Rightarrow g_{p,j} ({\bf u}) \leq 0.
$$

\noindent Note that an analogous result may be derived for the numerical constraints $g_j$, with the back-off of $\delta_e \kappa_{g,j}$ resulting in the same guarantee if we assume $g_j$ to be Lipschitz continuous with a Lipschitz constant of $\kappa_{g,j}$. However, more involved numerical techniques are also possible in this case \citep{Bunin2016}.

The exact manner in which one may employ (\ref{eq:backoff}) to allow for perturbation while robustly satisfying the constraints would almost certainly be method-dependent. The general principles, however, are quite simple. Because one has the guarantee that satisfying the constraint with a specified back-off allows for one to perturb safely anywhere in a $\delta_e$-ball around that point, it follows that one can choose $\delta_e$ \emph{a priori} and to then only perturb around those experimental iterates that satisfy the corresponding back-off, since all such information-gathering experiments will satisfy the constraints. This is now illustrated for a specific case-study example.

\subsubsection{Case Study: A Modifier-Adaptation Algorithm for Minimizing the Steady-State Production Cost of a Gold Cyanidation Leaching Process}
\label{sec:cyan}

One algorithm that depends heavily on perturbations is the modifier-adaptation method \citep{Brdys2005,Gao2005,Marchetti2009a}, which solves Problem (\ref{eq:mainprob}) by iteratively solving a corrected model-based problem:

\vspace{-3mm}
\begin{equation}\label{eq:mainprobMA}
\begin{array}{rll}
{\bf u}_{k+1} := & & \vspace{2mm} \\ {\rm arg} \mathop {{\rm{minimize}}}\limits_{\bf{u}} & \phi_{\hat p} ({\bf{u}},{\boldsymbol \theta}) + {\boldsymbol \lambda}_{\phi,k}^T {\bf u} &  \\
{\rm{subject}}\hspace{1mm}{\rm{to}} & g_{\hat p,j}({\bf{u}},{\boldsymbol \theta}) + \epsilon_{j,k} + {\boldsymbol \lambda}_{j,k}^T ( {\bf u} - {\bf u}_k ) \leq 0, & j = 1,...,n_{g_p} \vspace{1mm} \\
 & g_{j}({\bf{u}}) \leq 0, & j = 1,...,n_{g} \vspace{1mm} \\
& u_i^L \leq u_i \leq u_i^U, & i = 1,...,n_u,
\end{array}
\end{equation}

\noindent with $\epsilon$ and ${\boldsymbol \lambda}$ denoting the ``modifiers'' used to correct the local zero- and first-order errors of the model, respectively:

\vspace{-2mm}
\begin{equation}\label{eq:modifiers}
\begin{array}{l}
\epsilon_{j,k} := \alpha \left[ g_{p,j} ({\bf u}_k) - g_{\hat p,j} ({\bf u}_k, {\boldsymbol \theta}) \right] + (1-\alpha)\epsilon_{j,k-1} \vspace{1mm} \\
{\boldsymbol \lambda}_{j,k} := \alpha \left[ \nabla g_{p,j} ({\bf u}_k) - \nabla g_{\hat p,j} ({\bf u}_k, {\boldsymbol \theta}) \right] + (1-\alpha){\boldsymbol \lambda}_{j,k-1} \vspace{1mm} \\
{\boldsymbol \lambda}_{\phi,k} := \alpha \left[ \nabla \phi_{p} ({\bf u}_k) - \nabla \phi_{\hat p} ({\bf u}_k, {\boldsymbol \theta}) \right] + (1-\alpha){\boldsymbol \lambda}_{\phi,k-1}.
\end{array}
\end{equation}

\noindent It should be clear that this $\Gamma$ is a generalization of the constraint-adaptation algorithm employed in Section \ref{sec:poly}, in that we are now correcting the model \emph{derivatives} in addition to the constraint function values. As with constraint adaptation, the uncertain parameters of the model are not updated between iterations and are kept at some nominal, constant values, although a version of the algorithm that estimates and updates the parameters as well may certainly be used \citep{Brdys2005}.

As in the constraint-adaptation case, the zero-order modifier values ${\epsilon_{j,k}}$ may be easily obtained by physically measuring the experimental function and numerically evaluating the model at ${\bf u}_k$. The first-order corrections, however, require experimental derivative estimates, which may be obtained using a number of approaches \citep{Bunin2013a}. Here, we will employ the simple method of taking the difference quotients

\begin{subequations}\label{eq:findiff}
\begin{equation}\label{eq:findiffcost}
\frac{\partial f_{p}}{\partial u_i} \Big |_{{\bf u}_k} \approx \frac{f_{p}({\bf u}_k) - f_{p}({\bf u}_k - \delta_e {\bf e}_i)}{\delta_e},
\end{equation}
\noindent or, if $u_{k,i} - \delta_e < u_i^L$,
\begin{equation}\label{eq:findiffcon}
\frac{\partial f_{p}}{\partial u_i} \Big |_{{\bf u}_k} \approx \frac{f_{p}({\bf u}_k + \delta_e {\bf e}_i) - f_{p}({\bf u}_k)}{\delta_e},
\end{equation}
\end{subequations}

\noindent where ${\bf e}_i$ is the unit vector with unity at the $i^{\rm th}$ spot and $\delta_e$ is some non-zero value. The choice to use $\delta_e$ as the perturbation step is, of course, not coincidental and is intended to be used in synergy with (\ref{eq:backoff}) so as to ensure that all perturbations performed for estimating the derivatives satisfy the constraints.

The outline of the employed algorithm, where one starts at the iteration $k := 0$ with some ${\bf u}_0$ that satisfies the constraints on $g_{p,j}$ and $g_j$ with the additional back-offs of $\delta_e \kappa_{g_p,j}$ and $\delta_e \kappa_{g,j}$, is as follows:

\begin{enumerate}
\item Define the modifiers as in (\ref{eq:modifiers}), with the experimental gradients obtained by carrying out an additional $n_u$ experiments and then applying (\ref{eq:findiff}). 
\item Compute ${\bf u}_{k+1}$ by solving an augmented version of (\ref{eq:mainprobMA}) that accounts for the back-offs:

\begin{equation}\label{eq:mainprobMAlip}
\hspace{-10mm}\begin{array}{rll}
{\bf u}_{k+1} := & & \vspace{2mm} \\ {\rm arg} \mathop {{\rm{minimize}}}\limits_{\bf{u}} & \phi_{\hat p} ({\bf{u}},{\boldsymbol \theta}) + {\boldsymbol \lambda}_{\phi,k}^T {\bf u} &  \\
{\rm{subject}}\hspace{1mm}{\rm{to}} & g_{\hat p,j}({\bf{u}},{\boldsymbol \theta}) + \epsilon_{j,k} + {\boldsymbol \lambda}_{j,k}^T ( {\bf u} - {\bf u}_k ) \leq 0, & j = 1,...,n_{g_p} \vspace{1mm} \\
 & g_{j}({\bf{u}}) + \delta_e \kappa_{g,j} \leq 0, & j = 1,...,n_{g} \vspace{1mm} \\
& u_i^L \leq u_i \leq u_i^U, & i = 1,...,n_u \vspace{1mm} \\
& g_{p,j} ({\bf u}_k) + \kappa_{g_p,j}  \| {\bf u} - {\bf u}_k \|_2 + \delta_e \kappa_{g_p,j} \leq 0, & j = 1,...,n_{g_p},
\end{array}
\end{equation}

\item Apply the new decision variables to the experimental system, obtain the new measurements, augment the value of $k$, and return to Step 1.
\end{enumerate} 

Note that the modifications introduced in (\ref{eq:mainprobMAlip}) serve to ensure that the constraint function values -- experimental and numerical -- meet the back-offs required to ensure that the perturbations in the ball $\mathcal{B}_e$ around the future ${\bf u}_{k+1}$ satisfy the constraints.

This algorithm is applied to Problem P11 from the ExpOpt database \citep{ExpOpt}. The problem itself is adapted from the work of \cite{Jun2015} and deals with the minimization of the steady-state production cost of a gold cyanidation leaching process. Only the skeleton of the problem is provided here, with interested readers once more referred to the aforementioned references for more in-depth descriptions.

In our standard form, the problem may be written as

\begin{equation}\label{eq:cyanprob}
\begin{array}{rl}
\mathop {{\rm{minimize}}}\limits_{\bf{u}} & \phi_p ({\bf{u}}) := \phi_{\rm econ} ({\bf u}) \\
{\rm{subject}}\hspace{1mm}{\rm{to}} & \displaystyle g_{p,1}({\bf{u}}) := \frac{C_s ({\bf u}) - C_{s,0}}{C_{s,0}} + 0.75 \leq 0 \vspace{1mm} \\
& 10 \leq u_1 \leq 80 \vspace{1mm} \\
& 5 \leq u_2 \leq 20,
\end{array}
\end{equation}

\noindent with ${\bf u}$ consisting of two variables -- the flow rate of the sodium cyanide (in kilograms per hour) and the concentration of the dissolved oxygen in the liquid (in milligrams per kilogram). The function $\phi_{\rm econ}$ is the economical cost function to be minimized (in Chinese RMB per hour), while $C_s$ denotes the gold concentration in the ore (in milligrams per kilograms; $C_{s,0}$ denotes its initial value). The single experimental constraint ensures that the gold recovery is at least 75\% -- the solution, although not on this constraint, is nevertheless very close to it. The model used for the case study is that provided in the database and differs from the simulated ``real'' process by virtue of errors in the kinetic parameters. As in the previous case-study example, the decision variables are scaled down to a unit box, with the cost and constraint values divided by the scaling factors of 200 and 0.04, respectively. A $\delta_e$ value of 0.05 (for the scaled variables) is chosen, with the Lipschitz constant of the scaled problem set as $\kappa_{g_p,1} := 4$. The process starts operating at the suboptimal point of ${\bf u}_0 := (52, 18)$, which satisfies the experimental constraint with a slack greater than $\delta_e \kappa_{g_p,1}$. It is chosen not to filter the modifiers, with $\alpha := 1$ in (\ref{eq:modifiers}).

\begin{figure}
\begin{center}
\includegraphics[width=12cm]{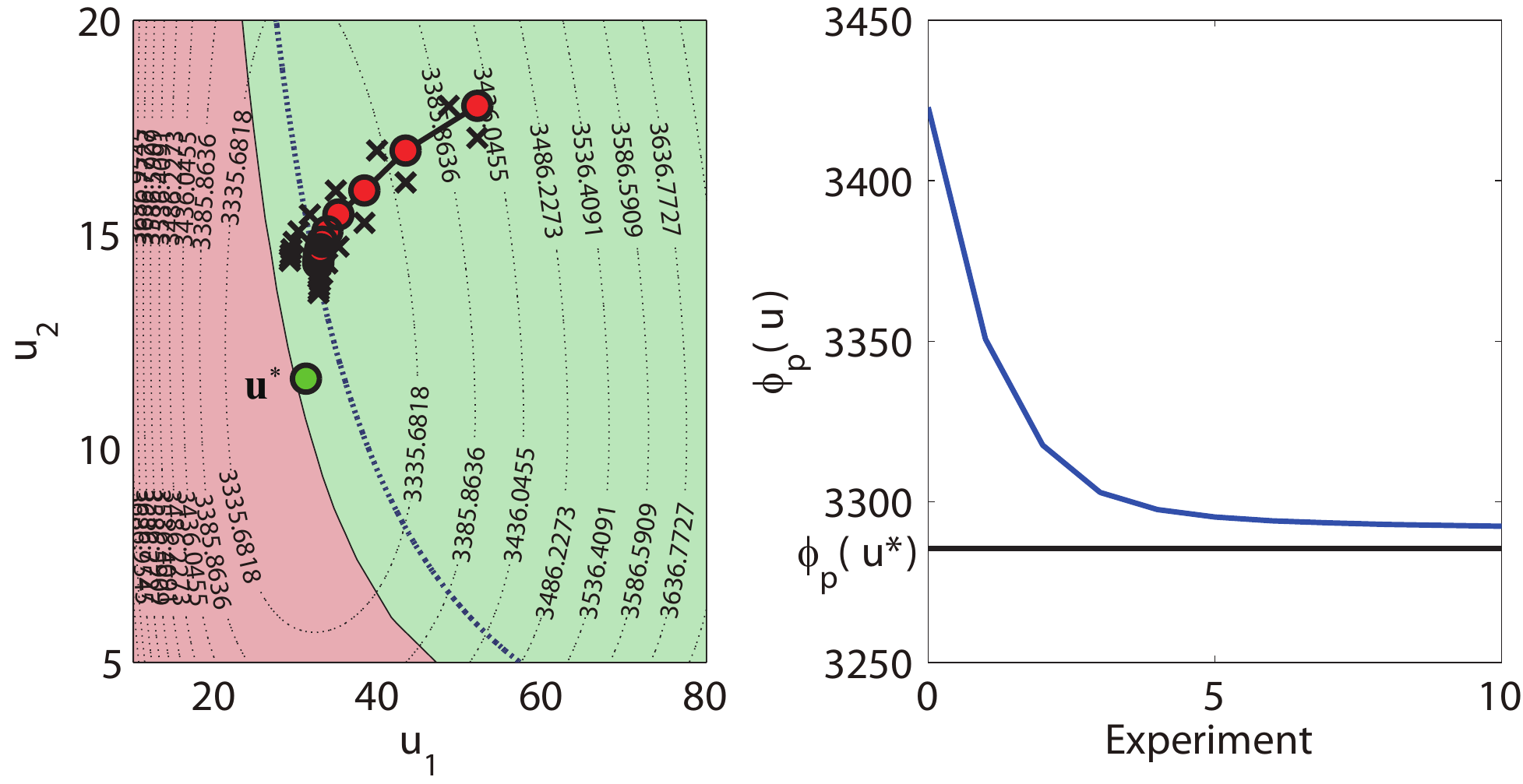}
\caption{Results obtained from applying the modifier-adaptation algorithm to Problem (\ref{eq:cyanprob}). The dashed line in the plot on the left shows the $\delta_e \kappa_{g_p,1} $ back-off sufficient to guarantee that the additional perturbations for derivative estimation satisfy the constraint. The cross marks show the perturbations, the data for which are omitted from the plot on the right.}
\label{fig:MAex}
\end{center}
\end{figure}

The results of applying the above algorithm to this problem are given in Fig. \ref{fig:MAex}. As expected, one sees that the main experimental iterates all satisfy the experimental constraint with the back-off, while the $\delta_e$-perturbations for derivative estimation all satisfy the original constraint. In this case, we see that the modifier-adaptation algorithm converges to the neighborhood of the optimum fairly quickly.

\subsection{Reducing the Effect of Uncertainty in the Experimental Function Values}
\label{sec:filter}

The Lipschitz constants may play an important supplementary role in the very realistic scenario where the experimental function values $f_p ({\bf u}_{\bar k}), \; \bar k = 0,...,k$ are subject to uncertainty, which may be due to (a) measurement noise or (b) the experimental quantity not being directly measured but estimated. In such cases, it is typical to work with a statistical interval in which the value $f_p ({\bf u}_{\bar k})$ is expected to lie. To obtain such an interval, one needs to make certain assumptions on the statistical nature of the uncertainty.

A common uncertainty model, and the one that will be considered here, is the additive

\begin{equation}\label{eq:meas}
\hat f_p^{\bar k} = f_p ({\bf u}_{\bar k}) + w_{\bar k},
\end{equation}

\noindent where $\hat f_p$ is the observed or estimated value and $w$ is the stochastic noise or error. By bounding the stochastic element in the probability sense\footnote{As a compromise between clarity and convenience, we will use the symbol $\mathop \leq \limits^P$ to mean ``less than or equal to \emph{with sufficiently high probability}''.}:

$$
\underline w_{\bar k} \mathop \leq \limits^P w_{\bar k} \mathop \leq \limits^P \overline w_{\bar k},
$$ 

\noindent it then becomes possible to compute a high-probability interval for $f_p ({\bf u}_{\bar k})$ by applying these bounds to a rearranged version of (\ref{eq:meas}):

\begin{equation}\label{eq:lipgenrob}
\underline f_p^{\bar k} := \hat f_p^{\bar k} - \overline w_{\bar k} \mathop \leq \limits^P f_p ({\bf u}_{\bar k}) \mathop \leq \limits^P \hat f_p^{\bar k} - \underline w_{\bar k} := \overline f_p^{\bar k}.
\end{equation}

Most importantly, these bounds allow us to obtain versions of (\ref{eq:congamma}) and (\ref{eq:backoff}) that are robust and do not require the unavailable exact function values:

\begin{equation}\label{eq:congammarob}
\overline g_{p,j}^k + \kappa_{g_p,j} \| {\bf u}_{k+1} - {\bf u}_k \|_2 \leq 0 \Rightarrow g_{p,j} ({\bf u}_{k+1}) \mathop \leq \limits^P 0,
\end{equation}

\begin{equation}\label{eq:backoffrob}
\overline g_{p,j}^k + \delta_e \kappa_{g_p,j}  \leq 0 \Rightarrow g_{p,j} ({\bf u}) \mathop \leq \limits^P 0, \;\; \forall {\bf u} \in \mathcal{B}_e \cap \mathcal{I},
\end{equation}

\noindent thereby making it possible to maintain the applicability of the methods in Sections \ref{sec:consat} and \ref{sec:excite} for the realistic case with uncertainty.

However, while these bounds may be fairly rigorous and useful, they can also prove too conservative when the variance of $w$ is large. A standard procedure for tightening the interval for a given ${\bf u}_{\bar k}$ is to resample the function value at ${\bf u}_{\bar k}$ multiple times and to take the mean of the values. Under reasonable assumptions on $w$, this then yields an estimator whose variance shrinks to 0 with an increasing number of samples \citep[\S 1.3]{Kay1993}. Alternatively, one may be able to obtain tighter bounding values by assuming a locally valid structure for $f_p$, such as linear or quadratic, and to fit this structure to the data via maximum-likelihood estimation. Provided that the structure chosen is correct, this too can be shown to yield estimators of $f_p$ values that go to the true value with a variance that decreases with an increasing amount of data \citep[\S 7]{Kay1993}. 

The alternative proposed here is to ``trim off'' some of the uncertainty when the resulting interval is actually larger than what would be allowed by the worst-case (Lipschitz-bound) changes in the function values from one decision-variable set to another. Consider first the robust version of (\ref{eq:lipgenLU}) for some $\tilde k \in \{ 0,...,k \}$:

$$
\underline f_p^{\tilde k} - \kappa \| {\bf u}_{\bar k} - {\bf u}_{\tilde k} \|_2 \mathop \leq \limits^P f_p ({\bf u}_{\bar k}) \mathop \leq \limits^P \overline f_p^{\tilde k} + \kappa \| {\bf u}_{\bar k} - {\bf u}_{\tilde k} \|_2,  
$$

\noindent of which (\ref{eq:lipgenrob}) are the particular case corresponding to $\tilde k \rightarrow \bar k$. We may, however, be able to improve the bounds by considering the tightest values over all $\{ 0,...,k \}$:

\begin{equation}\label{eq:lipgenrobL}
\underline f_p^{\bar k} := \mathop {\max} \limits_{{\tilde k} = 0,...,k} \left( \underline f_p^{\tilde k} - \kappa \| {\bf u}_{\bar k} - {\bf u}_{\tilde k} \|_2 \right) \mathop \leq \limits^P f_p ({\bf u}_{\bar k}),  
\end{equation}

\begin{equation}\label{eq:lipgenrobU}
f_p ({\bf u}_{\bar k}) \mathop \leq \limits^P \mathop {\min} \limits_{{\tilde k} = 0,...,k} \left( \overline f_p^{\tilde k} +  \kappa \| {\bf u}_{\bar k} - {\bf u}_{\tilde k} \|_2 \right) := \overline f_p^{\bar k}.
\end{equation}

After using (\ref{eq:lipgenrobL})-(\ref{eq:lipgenrobU}) to assign lower/upper bounding values to the measurements at ${\bf u}_0, {\bf u}_1,\hdots,{\bf u}_k$, one may run through the process again, since changes in some of the values may lead to additional changes in the others. The iterative scheme may be outlined as follows:

\begin{enumerate}
\item Set the nominal $\underline f_p^{\bar k}$ and $\overline f_p^{\bar k}$ for $\bar k = 0,...,k$ as given in (\ref{eq:lipgenrob}).
\item Refine these values by applying (\ref{eq:lipgenrobL})-(\ref{eq:lipgenrobU}) for $\bar k = 0,...,k$.
\item If the largest refinement obtained is negligible (e.g., $< 10^{-6}$), terminate. Otherwise, return to Step 2.
\end{enumerate}

\begin{figure}
\begin{center}
\includegraphics[width=12cm]{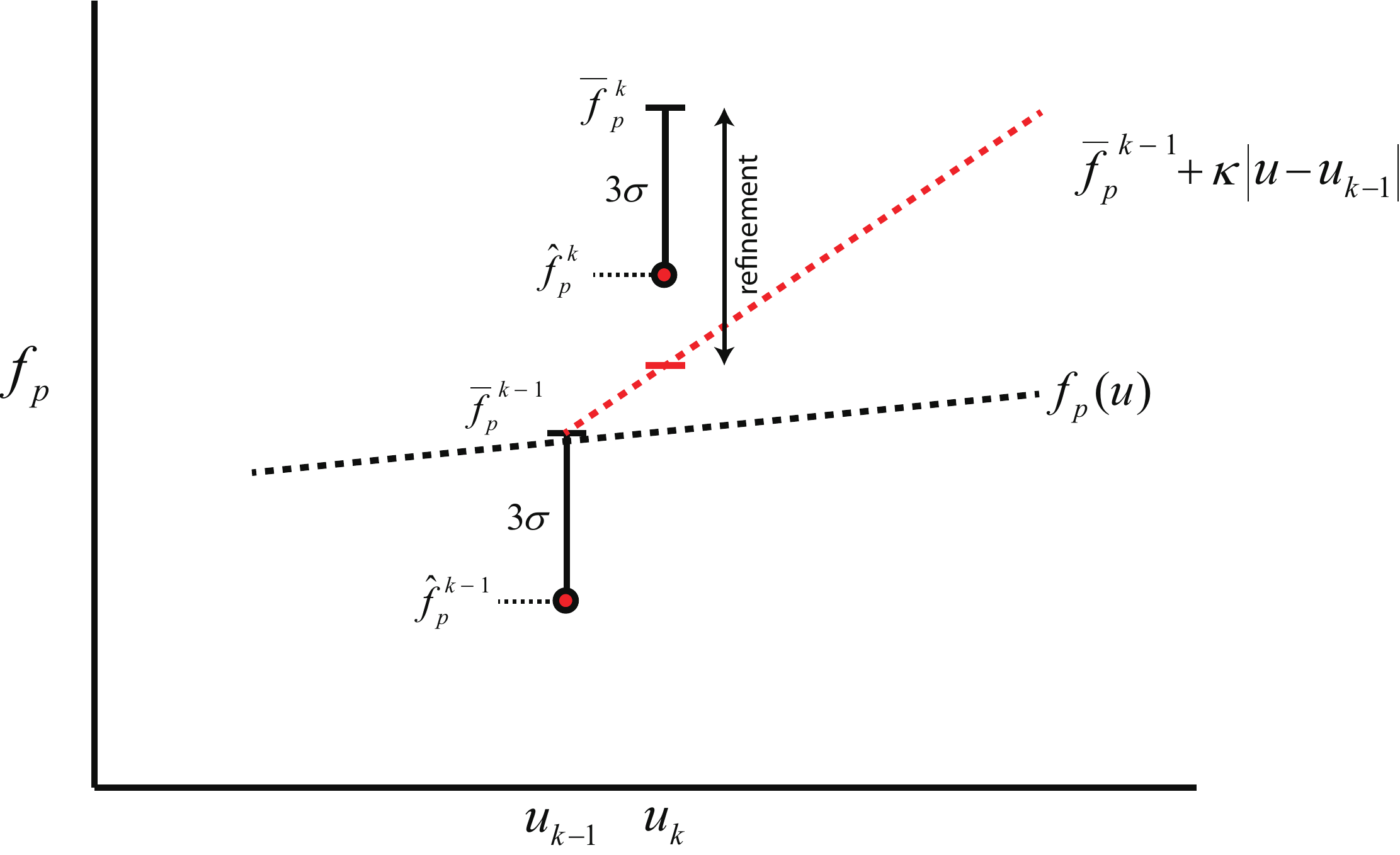}
\caption{Illustration of how Lipschitz bounds may be used to tighten uncertainty bounds on the true function values. The red dashed line shows the Lipschitz bound.}
\label{fig:liprefine}
\end{center}
\end{figure}

A conceptual illustration of reducing the uncertainty via this procedure is given in Fig. \ref{fig:liprefine}, where we consider a simple one-dimensional case with white Gaussian noise of variance $\sigma^2$, and high-probability bounds of $-3 \sigma$ and $3 \sigma$ on the noise elements. As shown in the figure, it happens that a very poor nominal upper bounding value is obtained at $u_k$ after a very tight upper bounding value has been obtained at the previous experiment of $u_{k-1}$. Because of the Lipschitz bound, the tightness at the latter may then be partially ``inherited'' by the former, thereby leading to significant refinement in the upper bounding value at $u_k$. Naturally, we may expect this method to lose its effectiveness when the decision-variable points are further apart and when the Lipschitz constants used are more conservative.

\subsubsection{Case Study: Gold Cyanidation Leaching Process with Noise}

Let us again consider the case-study example of Section \ref{sec:cyan}, this time supposing a more realistic scenario where all of the measurements obtained, for both the cost and constraint, are corrupted by noise elements from the PDF of $\mathcal{N} (0, \sigma^2)$. A value of $\sigma := 0.07$ is used here, and the additive corruption is applied to the problem post-scaling. The nominal lower and upper bounding values on the noise elements are then chosen as $\underline w := -3 \sigma$ and $\overline w := 3 \sigma$ for all measurements at all iterations, from which the nominal lower and upper bounding values on the true function values are obtained as in (\ref{eq:lipgenrob}).

Two realizations of the modifier-adaptation algorithm are compared here, with both still setting the modifiers as in (\ref{eq:modifiers})-(\ref{eq:findiff}), but with the $f_p$ function values replaced by their measured or estimated $\hat f_p$ analogues. To robustify the constraint-satisfaction guarantees against the uncertainty effects, the numerical optimization problem of (\ref{eq:mainprobMAlip}) is modified to

\begin{equation}\label{eq:mainprobMAlip2}
\hspace{-10mm}\begin{array}{rll}
{\bf u}_{k+1} := & & \vspace{2mm} \\ {\rm arg} \mathop {{\rm{minimize}}}\limits_{\bf{u}} & \phi_{\hat p} ({\bf{u}},{\boldsymbol \theta}) + {\boldsymbol \lambda}_{\phi,k}^T {\bf u} &  \\
{\rm{subject}}\hspace{1mm}{\rm{to}} & g_{\hat p,j}({\bf{u}},{\boldsymbol \theta}) + \epsilon_{j,k} + {\boldsymbol \lambda}_{j,k}^T ( {\bf u} - {\bf u}_k ) \leq 0, & j = 1,...,n_{g_p} \vspace{1mm} \\
 & g_{j}({\bf{u}}) + \delta_e \kappa_{g,j} \leq 0, & j = 1,...,n_{g} \vspace{1mm} \\
& u_i^L \leq u_i \leq u_i^U, & i = 1,...,n_u \vspace{1mm} \\
& \overline g_{p,j}^k + \kappa_{g_p,j}  \| {\bf u} - {\bf u}_k \|_2 + \delta_e \kappa_{g_p,j} \leq 0, & j = 1,...,n_{g_p}.
\end{array}
\end{equation}

In the first algorithm, the Lipschitz constants are \emph{not} used to refine the bounding values, with the $\overline g_{p,j}^k$ value obtained from (\ref{eq:lipgenrob}) used in (\ref{eq:mainprobMAlip2}). In the second algorithm, the Lipschitz-bound refinement procedure as outlined in the preceding section is run through all of the measurements -- the ``key'' iterations $0,...,k$ \emph{and} the additional perturbations around the key iterations used for derivative estimation. The resulting $\overline g_{p,j}^{k}$ is then used in (\ref{eq:mainprobMAlip2}). Additionally, the measured values $\hat \phi_p$ and $\hat g_{p,j}$ that are used in estimating the modifiers and derivatives are trimmed in the case that they do not satisfy the tightened bounding values (e.g., $\hat \phi^k_p \rightarrow \overline \phi^k_p$ if $\hat \phi^k_p > \overline \phi^k_p$). The Lipschitz constant for the cost is chosen as $\kappa_{\phi} := 5$.

So as to see the average effect that including the Lipschitz-based tightening has on the performance of the optimization algorithm, we run the algorithms in parallel for 1000 noise realizations and evaluate the average difference between the corresponding cost function values,

$$
\Delta \phi_{\rm ave} := \frac{1}{k_f+1} \sum_{\bar k = 0}^{k_f} (\phi_{p,\bar k}^A - \phi_{p,\bar k}^B),
$$

\noindent for each realization, with $\phi_p^A$ and $\phi_p^B$ denoting the cost function values obtained by the original (without Lipschitz tightening) and modified (with Lipschitz tightening) algorithms, respectively. The scalar $k_f$ is used to denote the number of main experiments run, excluding the perturbations and the initial experiment, and is set as 30 here. A positive $\Delta \phi_{\rm ave}$ value thus suggests an improvement over the original algorithm. 

\begin{figure}
\begin{center}
\includegraphics[width=12cm]{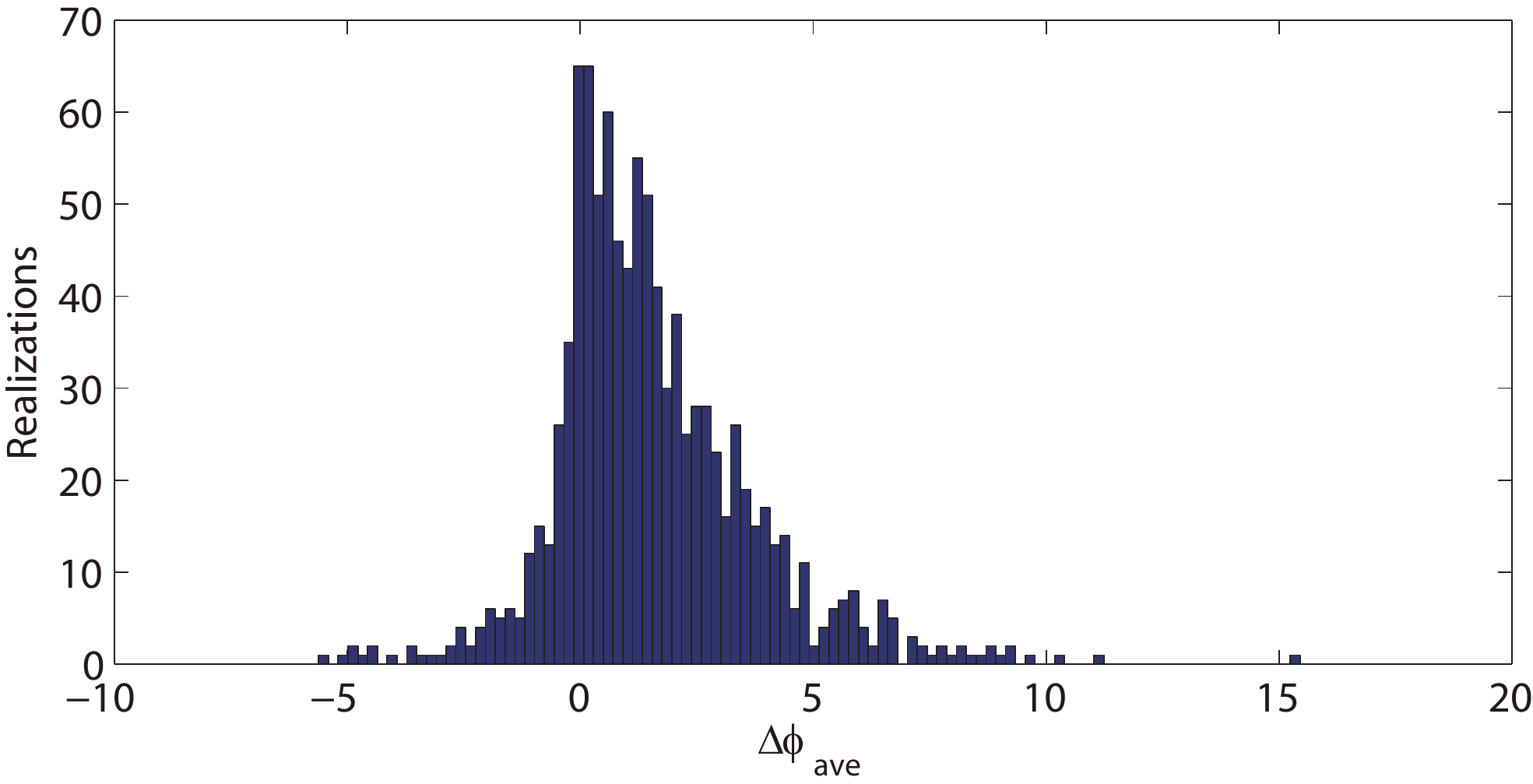}
\caption{Distribution of the average differences between the cost function values obtained by the two modifier-adaptation algorithms.}
\label{fig:hist}
\end{center}
\end{figure}

Plotting the results of the thousand trials in a histogram (Fig. \ref{fig:hist}), we notice a very positive trend, with noticeably better results obtained by the second algorithm.

\section{Improving the Lipschitz Bounds}
\label{sec:refine}

A significant price to pay for the use of the Lipschitz bounds of (\ref{eq:lipgenLU}) is their potential conservatism. Namely, if the constant $\kappa$ is chosen so that the bounds are satisfied with a very large margin, one is faced with several drawbacks:

\begin{itemize}
\item as already mentioned in Section \ref{sec:consat}, the constraint-satisfaction condition (\ref{eq:congamma}) may become too restrictive, forcing the algorithm to take overly small steps;
\item the back-off, $\delta_e  \kappa_{p,j}$, needed for safe perturbation in (\ref{eq:backoff}) may become too large and difficult to satisfy, and if satisfied may introduce a large degree of suboptimality;
\item the Lipschitz-based tightening of the bounding values $\underline f_p$ and $\overline f_p$ obtained in Section \ref{sec:filter} becomes negligible or disappears altogether.
\end{itemize}

Fortunately, though the function $f_p$ may be experimental and thus unknown, scattered bits of information about its behavior might be available and could be exploited to yield different Lipschitz bounds that are often much tighter.

In this section, we consider the following types of additional information:

\begin{itemize}
\item directional lower and upper Lipschitz constants,

$$
\underline \kappa_i \leq \frac{\partial f_p}{\partial u_i} \Big |_{\bf u} \leq \overline \kappa_i, \; \forall {\bf u} \in \mathcal{I};
$$

\item local Lipschitz constants over the subspace 

$$
\mathcal{I}_{{\bf u}_a}^{{\bf u}_b} = \{ {\bf u} : u_{a,i} \leq u_i \leq u_{b,i}, \; i = 1,...,n_u \} \subseteq \mathcal{I}_{{\bf u}^L}^{{\bf u}^U} \subseteq \mathcal{I},
$$

\noindent defined as

\begin{equation}\label{eq:liplocal}
\underline \kappa_i^{{\bf u}_a, {\bf u}_b} \leq \frac{\partial f_p}{\partial u_i} \Big |_{\bf u} \leq \overline \kappa_i^{{\bf u}_a, {\bf u}_b}, \; \forall {\bf u} \in \mathcal{I}_{{\bf u}_a}^{{\bf u}_b};
\end{equation}

\item the local convexity and concavity properties of $f_p$, where it is known that $f_p$ is locally convex or concave in certain variables over $\mathcal{I}_{{\bf u}_a}^{{\bf u}_b}$, with the indices of these variables denoted by the sets $I_{\rm cvx}^{{\bf u}_a, {\bf u}_b}$ and $I_{\rm ccv}^{{\bf u}_a, {\bf u}_b}$, respectively;

\item lower and upper bounds on the local derivatives of $f_p$:

\begin{equation}\label{eq:derbounds}
\underline \kappa_i^{{\bf u}_a, {\bf u}_b} \leq \nabla \underline f_{p,i}^{{\bf u}_a} \leq \frac{\partial f_p}{\partial u_i} \Big |_{{\bf u}_a} \leq \nabla \overline f_{p,i}^{{\bf u}_a} \leq \overline \kappa_i^{{\bf u}_a, {\bf u}_b}.
\end{equation}

\end{itemize}

By integrating the information above, one is able to account for scenarios that would be completely ignored when using (\ref{eq:lipgenLU}). For example, by using $2n_u$ Lipschitz constants instead of just one, it becomes possible to differentiate between those decision variables that have a large effect on $f_p$ and those whose effect is known to be negligible. Additionally, by having both lower ($\underline \kappa_i$) and upper ($\overline \kappa_i$) constants, one can differentiate between the sensitivities as being either positive or negative -- e.g., if it is known that increasing a certain $u_i$ will \emph{always} increase the value of $f_p$, the corresponding $\underline \kappa_i$ can be set as 0, which, as will be seen, may tighten the corresponding Lipschitz bound significantly.  

By limiting the analysis to a subspace $\mathcal{I}_{{\bf u}_a}^{{\bf u}_b} \subseteq \mathcal{I}$, it becomes possible to accomodate those problems where the sensitivity of $f_p$ to some $u_i$ is known to change drastically over $\mathcal{I}$. Of particular concern may be those $f_p$ that have an inversely proportional relationship with $u_i$, and whose values may explode as $u_i$ approaches 0. While it might make sense to include the 0-point in the definition of $\mathcal{I}$, actually operating there might be unlikely, and so it would be unnecessarily conservative to use the very large Lipschitz constants corresponding to the region around the 0-point when smaller, more reasonable values could be used instead. Considering the local subspace $\mathcal{I}_{{\bf u}_a}^{{\bf u}_b}$ provides us with a convenient mathematical formulation for handling such cases.

Finally, the inclusion of convexity/concavity properties, when coupled with derivative bounds, gives greater flexibility for those problems where $f_p$ is known with certainty to behave in a convex/concave manner when only certain variables are changed (the other variables being kept the same). The relationship between the power and current in a fuel-cell system \citep{Marchetti2009}, known to be concave, is one example that the authors have encountered. Mathematically, exploiting such relationships allows us to forgo using Lipschitz constants for certain directions and to use the less conservative derivative bounds instead. Some results for obtaining such bounds for an experimental function $f_p$ have been reported in \cite{Bunin2013a}.

Combining the information above, let us now derive an alternative to (\ref{eq:lipgenLU}).

\begin{thm}[Alternate Lipschitz Bounds]
Let the experimental function $f_p$ be continuously differentiable ($\mathcal{C}^1$) over an open set containing $\mathcal{I}_{{\bf u}_a}^{{\bf u}_b}$, with its derivatives bounded as in (\ref{eq:liplocal})-(\ref{eq:derbounds}). The bounds

\begin{equation}\label{eq:lipgenL}
\hspace{-15mm} f_p ({\bf u}_a) + \hspace{-3mm} \sum_{i \in I_{\rm cvx}^{{\bf u}_a, {\bf u}_b}} \hspace{-3mm} \mathop {\min} \left[ \begin{array}{l} \nabla \underline f_{p,i}^{{\bf u}_a} (u_{b,i} - u_{a,i}), \vspace{1mm} \\  \nabla \overline f_{p,i}^{{\bf u}_a} (u_{b,i} - u_{a,i}) \end{array} \right] + \hspace{-3mm} \sum_{i \not \in I_{\rm cvx}^{{\bf u}_a, {\bf u}_b}} \hspace{-3mm} \mathop {\min} \left[ \begin{array}{l} \underline \kappa_{i}^{{\bf u}_a,{\bf u}_b} (u_{b,i} - u_{a,i}), \vspace{1mm} \\  \overline \kappa_{i}^{{\bf u}_a,{\bf u}_b} (u_{b,i} - u_{a,i}) \end{array} \right] \leq f_p ({\bf u}_b),
\end{equation}

\begin{equation}\label{eq:lipgenU}
\hspace{-15mm} f_p({\bf u}_b) \leq f_p ({\bf u}_a) + \hspace{-3mm} \sum_{i \in I_{\rm ccv}^{{\bf u}_a, {\bf u}_b}} \hspace{-3mm} \mathop {\max} \left[ \begin{array}{l} \nabla \underline f_{p,i}^{{\bf u}_a} (u_{b,i} - u_{a,i}), \vspace{1mm} \\  \nabla \overline f_{p,i}^{{\bf u}_a} (u_{b,i} - u_{a,i}) \end{array} \right] + \hspace{-3mm} \sum_{i \not \in I_{\rm ccv}^{{\bf u}_a, {\bf u}_b}} \hspace{-3mm} \mathop {\max} \left[ \begin{array}{l} \underline \kappa_{i}^{{\bf u}_a,{\bf u}_b} (u_{b,i} - u_{a,i}), \vspace{1mm} \\  \overline \kappa_{i}^{{\bf u}_a,{\bf u}_b} (u_{b,i} - u_{a,i}) \end{array} \right]
\end{equation}

\noindent are then satisfied for all ${\bf u}_a, {\bf u}_b \in \mathcal{I}_{{\bf u}_a}^{{\bf u}_b}$.

\end{thm}
\begin{pf}
We will only prove (\ref{eq:lipgenU}) and leave (\ref{eq:lipgenL}) to the reader, as the two proofs are symmetric. First, let us write the evolution of $f_p$ from ${\bf u}_a$ to ${\bf u}_b$ as occurring in ``two steps'', the first step in the variables in which $f_p$ is concave and the second step in the variables in which $f_p$ is not. For clarity, we will partition ${\bf u}$ as ${\bf u} = ({\bf v},{\bf z})$, where ${\bf v}$ are the variables $u_i, i \in I_{\rm ccv}^{{\bf u}_a, {\bf u}_b}$ and ${\bf z}$ are the variables $u_i, i \not \in I_{\rm ccv}^{{\bf u}_a, {\bf u}_b}$. The two-step evolution may be written as

$$
f_p ({\bf v}_b, {\bf z}_b) = f_p ({\bf v}_a, {\bf z}_a) + \left[ f_p ({\bf v}_b, {\bf z}_a) - f_p ({\bf v}_a, {\bf z}_a) \right] + \left[ f_p ({\bf v}_b, {\bf z}_b) - f_p ({\bf v}_b, {\bf z}_a) \right].
$$

The proof consists in upper-bounding the two differences in the brackets. For the first difference, we call upon the first-order condition of concavity \citep[\S 3]{Boyd2008}:

$$
f_p ({\bf v}_b, {\bf z}_a) - f_p ({\bf v}_a, {\bf z}_a) \leq \sum_{i \in I_{\rm ccv}^{{\bf u}_a, {\bf u}_b}} \frac{\partial f_p}{\partial u_i} \Big |_{{\bf u}_a} (u_{b,i} - u_{a,i}).
$$

To further bound the right-hand side, we then consider two cases,

\begin{equation}\label{eq:posneg}
\begin{array}{l}
\displaystyle u_{b,i} - u_{a,i} \geq 0 \Leftrightarrow \frac{\partial f_p}{\partial u_i} \Big |_{{\bf u}_a} (u_{b,i} - u_{a,i})  \leq \nabla \overline f_{p,i}^{{\bf u}_a} (u_{b,i} - u_{a,i}) \vspace{1mm} \\
\displaystyle u_{b,i} - u_{a,i} \leq 0 \Leftrightarrow \frac{\partial f_p}{\partial u_i} \Big |_{{\bf u}_a} (u_{b,i} - u_{a,i})  \leq \nabla \underline f_{p,i}^{{\bf u}_a} (u_{b,i} - u_{a,i}),
\end{array}
\end{equation}

\noindent and account for the worst case via the maximum operator to obtain

\begin{equation}\label{eq:firstbound}
f_p ({\bf v}_b, {\bf z}_a) - f_p ({\bf v}_a, {\bf z}_a) \leq \sum_{i \in I_{\rm ccv}^{{\bf u}_a, {\bf u}_b}} \hspace{-3mm} \mathop {\max} \left[ \begin{array}{l} \nabla \underline f_{p,i}^{{\bf u}_a} (u_{b,i} - u_{a,i}), \vspace{1mm} \\  \nabla \overline f_{p,i}^{{\bf u}_a} (u_{b,i} - u_{a,i}) \end{array} \right].
\end{equation}

To place a bound on $f_p ({\bf v}_b, {\bf z}_b) - f_p ({\bf v}_b, {\bf z}_a)$, we limit our analysis to a single dimension and then apply the mean-value theorem. First, consider the function $f_p$ over the line connecting $({\bf v}_b, {\bf z}_a)$ and $({\bf v}_b, {\bf z}_b)$, parameterized as

$$
\tilde f_p ({\bf v}_b, \gamma) = f_p ({\bf v}_b, {\bf z} (\gamma)),
$$ 

\noindent where

$$
{\bf z} (\gamma) = {\bf z}_a + \gamma ({\bf z}_b - {\bf z}_a), \; \gamma \in [0,1].
$$

Because $f_p$ is $\mathcal{C}^1$, $\tilde f_p$ is as well, thereby allowing us to apply the Taylor series expansion between $\gamma = 0$ and $\gamma = 1$, together with the mean-value theorem \citep[\S 4.7]{Korn2000}, to state that

$$
\exists \tilde \gamma \in (0,1) : \tilde f_p ({\bf v}_b, 1) = \tilde f_p ({\bf v}_b,0) + \frac{\partial \tilde f_p}{\partial \gamma} \Big |_{{\bf v}_b, \tilde \gamma}.
$$

The derivative term may then be defined in terms of the original function $f_p$ by applying the chain rule:

$$
\frac{\partial \tilde f_p}{\partial \gamma} \Big |_{{\bf v}_b, \gamma} = \sum_{i \not \in I_{\rm ccv}^{{\bf u}_a, {\bf u}_b}} \frac{\partial f_p}{\partial u_i} \Big |_{{\bf v}_b, {\bf z}(\gamma)} \frac{d u_i}{d \gamma} \Big |_{\gamma} = \sum_{i \not \in I_{\rm ccv}^{{\bf u}_a, {\bf u}_b}} \frac{\partial f_p}{\partial u_i} \Big |_{{\bf v}_b, {\bf z}(\gamma)} (u_{b,i} - u_{a,i}).
$$

Noting that $\tilde f_p ({\bf v}_b,0) = f_p ({\bf v}_b, {\bf z}_a)$ and $\tilde f_p ({\bf v}_b,1) = f_p ({\bf v}_b, {\bf z}_b)$, we thus obtain that

$$
\exists \tilde \gamma \in (0,1) : f_p ({\bf v}_b, {\bf z}_b) = f_p ({\bf v}_b, {\bf z}_a) + \sum_{i \not \in I_{\rm ccv}^{{\bf u}_a, {\bf u}_b}} \frac{\partial f_p}{\partial u_i} \Big |_{{\bf v}_b, {\bf z}(\tilde \gamma)} (u_{b,i} - u_{a,i}).
$$

Because $({\bf v}_b, {\bf z}(\tilde \gamma)) \in \mathcal{I}_{{\bf u}_a}^{{\bf u}_b}$, the constants of (\ref{eq:liplocal}) may be used to bound the derivatives and, using the same logic as in (\ref{eq:posneg}), we go on to obtain

$$
\frac{\partial f_p}{\partial u_i} \Big |_{{\bf v}_b, {\bf z}(\tilde \gamma)} (u_{b,i} - u_{a,i}) \leq \mathop {\max} \left[ \begin{array}{l} \underline \kappa_{i}^{{\bf u}_a,{\bf u}_b} (u_{b,i} - u_{a,i}), \vspace{1mm} \\  \overline \kappa_{i}^{{\bf u}_a,{\bf u}_b} (u_{b,i} - u_{a,i}) \end{array} \right],
$$

\noindent which then leads to

\begin{equation}\label{eq:secondbound}
f_p ({\bf v}_b, {\bf z}_b) - f_p ({\bf v}_b, {\bf z}_a) \leq \sum_{i \not \in I_{\rm ccv}^{{\bf u}_a, {\bf u}_b}} \hspace{-3mm} \mathop {\max} \left[ \begin{array}{l} \underline \kappa_{i}^{{\bf u}_a,{\bf u}_b} (u_{b,i} - u_{a,i}), \vspace{1mm} \\  \overline \kappa_{i}^{{\bf u}_a,{\bf u}_b} (u_{b,i} - u_{a,i}) \end{array} \right].
\end{equation}

Combining (\ref{eq:firstbound}) and (\ref{eq:secondbound}) yields the desired result. \qed

\end{pf}

Next, we consider the analogue to (\ref{eq:backoff}) that follows from (\ref{eq:lipgenU}).

\begin{corol}[Alternate Back-off for Local Perturbations]\label{corol:backoff}
Let $g_{p,j}$ be $\mathcal{C}^1$ over an open set containing $\mathcal{I}_{{\bf u}_k-{\boldsymbol \delta}_e}^{{\bf u}_k + {\boldsymbol \delta}_e}$, where ${\boldsymbol \delta}_e \in \mathbb{R}^{n_u}$ is the vector of $\delta_e$ values, and let the derivatives be bounded as in (\ref{eq:liplocal})-(\ref{eq:derbounds}) with $\mathcal{I}_{{\bf u}_a}^{{\bf u}_b} \rightarrow \mathcal{I}_{{\bf u}_k-{\boldsymbol \delta}_e}^{{\bf u}_k + {\boldsymbol \delta}_e}$ and $f_p \rightarrow g_{p,j}$. Defining the mixed vector of greatest absolute Lipschitz constants and derivative bounds as

$$
\tilde \kappa_{p,ji} = \left\{ \begin{array}{ll} \displaystyle \mathop {\max} \left( \Big | \nabla \underline g_{p,ji}^{{\bf u}_k} \Big |, \Big | \nabla \overline g_{p,ji}^{{\bf u}_k} \Big | \right), & i \in I_{{\rm ccv},j}^{{\bf u}_k-{\boldsymbol \delta}_e,{\bf u}_k+{\boldsymbol \delta}_e}, \vspace{1mm} \\ \displaystyle \mathop {\max} \left( \Big |\underline \kappa_{p,ji}^{{\bf u}_k-{\boldsymbol \delta}_e,{\bf u}_k+{\boldsymbol \delta}_e} \Big |, \Big |\overline \kappa_{p,ji}^{{\bf u}_k-{\boldsymbol \delta}_e,{\bf u}_k+{\boldsymbol \delta}_e} \Big | \right), & i \not \in I_{{\rm ccv},j}^{{\bf u}_k-{\boldsymbol \delta}_e,{\bf u}_k+{\boldsymbol \delta}_e}, \end{array} \right .
$$

\noindent it follows that

\begin{equation}\label{eq:backoff2}
g_{p,j} ({\bf u}_k) + \delta_e \| \tilde {\boldsymbol \kappa}_{p,j} \|_2 \leq 0 \Rightarrow g_{p,j} ({\bf u}) \leq 0,  \; \forall {\bf u} \in \mathcal{B}_e.
\end{equation}

\end{corol}
\begin{pf}
Applying (\ref{eq:lipgenU}), it follows that

\begin{equation}\label{eq:lipbounde}
\begin{array}{l}
\displaystyle g_{p,j} ({\bf u}) \leq g_{p,j} ({\bf u}_k) + \sum_{i \in I_{{\rm ccv},j}^{{\bf u}_k-{\boldsymbol \delta}_e,{\bf u}_k+{\boldsymbol \delta}_e}} \hspace{-3mm} \mathop {\max} \left[ \begin{array}{l} \nabla \underline g_{p,ji}^{{\bf u}_k} (u_{i} - u_{k,i}), \vspace{1mm} \\  \nabla \overline g_{p,ji}^{{\bf u}_k} (u_{i} - u_{k,i}) \end{array} \right] \vspace{1mm} \\ 
\hspace{35mm} \displaystyle + \hspace{-3mm} \sum_{i \not \in I_{{\rm ccv},j}^{{\bf u}_k-{\boldsymbol \delta}_e,{\bf u}_k+{\boldsymbol \delta}_e}} \hspace{-3mm} \mathop {\max} \left[ \begin{array}{l} \underline \kappa_{p,ji}^{{\bf u}_k-{\boldsymbol \delta}_e,{\bf u}_k+{\boldsymbol \delta}_e} (u_{i} - u_{k,i}), \vspace{1mm} \\  \overline \kappa_{p,ji}^{{\bf u}_k-{\boldsymbol \delta}_e,{\bf u}_k+{\boldsymbol \delta}_e} (u_{i} - u_{k,i}) \end{array} \right]
\end{array}
\end{equation}

\noindent holds for all ${\bf u} \in \mathcal{B}_e \subseteq \mathcal{I}_{{\bf u}_k-{\boldsymbol \delta}_e}^{{\bf u}_k + {\boldsymbol \delta}_e}$. Let us use this bound to define a \emph{Lipschitz polytope} centered at ${\bf u}_k$:

$$
\mathcal{L}_j = \left\{ {\bf u} : \begin{array}{l} \displaystyle g_{p,j} ({\bf u}_k) + \sum_{i \in I_{{\rm ccv},j}^{{\bf u}_k-{\boldsymbol \delta}_e,{\bf u}_k+{\boldsymbol \delta}_e}} \hspace{-3mm} \mathop {\max} \left[ \begin{array}{l} \nabla \underline g_{p,ji}^{{\bf u}_k} (u_{i} - u_{k,i}), \vspace{1mm} \\  \nabla \overline g_{p,ji}^{{\bf u}_k} (u_{i} - u_{k,i}) \end{array} \right] \vspace{1mm} \\ 
\hspace{15mm} \displaystyle + \hspace{-3mm} \sum_{i \not \in I_{{\rm ccv},j}^{{\bf u}_k-{\boldsymbol \delta}_e,{\bf u}_k+{\boldsymbol \delta}_e}} \hspace{-3mm} \mathop {\max} \left[ \begin{array}{l} \underline \kappa_{p,ji}^{{\bf u}_k-{\boldsymbol \delta}_e,{\bf u}_k+{\boldsymbol \delta}_e} (u_{i} - u_{k,i}), \vspace{1mm} \\  \overline \kappa_{p,ji}^{{\bf u}_k-{\boldsymbol \delta}_e,{\bf u}_k+{\boldsymbol \delta}_e} (u_{i} - u_{k,i}) \end{array} \right] \leq 0  \end{array} \right\}.
$$

The remainder of the proof is essentially a ``busier'' version of the proof of Theorem 1 in \cite{Bunin2016}, with the same concepts at work. From (\ref{eq:lipbounde}), it should be clear that any ${\bf u}$ belonging to both $\mathcal{B}_e$ and $\mathcal{L}_j$ will satisfy $g_{p,j} ({\bf u}) \leq 0$, and so we proceed to identify the back-off on $g_{p,j} ({\bf u}_k)$ that would enforce $\mathcal{B}_e \subseteq \mathcal{L}_j$, and thus ${\bf u} \in \mathcal{B}_e \Rightarrow {\bf u} \in \mathcal{L}_j$ (i.e., by inscribing $\mathcal{B}_e$ in $\mathcal{L}_j$). This is done analytically by building progressively smaller subsets of $\mathcal{L}_j$, all of which include $\mathcal{B}_e$. 

From the algebraic relation $\mathop {\max} (\underline x y, \overline x y) \leq \mathop {\max} (| \underline x |, | \overline x |) |y|$, we have

$$
\mathop {\max} \left[ \begin{array}{l} \nabla \underline g_{p,ji}^{{\bf u}_k} (u_{i} - u_{k,i}), \vspace{1mm} \\  \nabla \overline g_{p,ji}^{{\bf u}_k} (u_{i} - u_{k,i}) \end{array} \right] \leq \mathop {\max} \left( \Big | \nabla \underline g_{p,ji}^{{\bf u}_k} \Big |, \Big | \nabla \overline g_{p,ji}^{{\bf u}_k} \Big | \right) | u_i - u_{k,i} |,
$$

$$
\mathop {\max} \left[ \begin{array}{l} \underline \kappa_{p,ji}^{{\bf u}_k-{\boldsymbol \delta}_e,{\bf u}_k+{\boldsymbol \delta}_e} (u_{i} - u_{k,i}), \vspace{1mm} \\  \overline \kappa_{p,ji}^{{\bf u}_k-{\boldsymbol \delta}_e,{\bf u}_k+{\boldsymbol \delta}_e} (u_{i} - u_{k,i}) \end{array} \right] \leq \mathop {\max} \left( \Big |\underline \kappa_{p,ji}^{{\bf u}_k-{\boldsymbol \delta}_e,{\bf u}_k+{\boldsymbol \delta}_e} \Big |, \Big |\overline \kappa_{p,ji}^{{\bf u}_k-{\boldsymbol \delta}_e,{\bf u}_k+{\boldsymbol \delta}_e} \Big | \right) | u_i - u_{k,i} |,
$$

\noindent whence follows

\begin{equation}\label{eq:proofineq}
\hspace{-25mm}\begin{array}{l}
\displaystyle \sum_{i \in I_{{\rm ccv},j}^{{\bf u}_k-{\boldsymbol \delta}_e,{\bf u}_k+{\boldsymbol \delta}_e}} \hspace{-3mm} \mathop {\max} \left[ \begin{array}{l} \nabla \underline g_{p,ji}^{{\bf u}_k} (u_{i} - u_{k,i}), \vspace{1mm} \\  \nabla \overline g_{p,ji}^{{\bf u}_k} (u_{i} - u_{k,i}) \end{array} \right] + \sum_{i \not \in I_{{\rm ccv},j}^{{\bf u}_k-{\boldsymbol \delta}_e,{\bf u}_k+{\boldsymbol \delta}_e}} \hspace{-3mm} \mathop {\max} \left[ \begin{array}{l} \underline \kappa_{p,ji}^{{\bf u}_k-{\boldsymbol \delta}_e,{\bf u}_k+{\boldsymbol \delta}_e} (u_{i} - u_{k,i}), \vspace{1mm} \\  \overline \kappa_{p,ji}^{{\bf u}_k-{\boldsymbol \delta}_e,{\bf u}_k+{\boldsymbol \delta}_e} (u_{i} - u_{k,i}) \end{array} \right] \vspace{1mm} \\
\hspace{75mm}\displaystyle \leq \sum_{i = 1}^{n_u} \tilde \kappa_{p,ji} | u_i - u_{k,i} | = \tilde {\boldsymbol \kappa}_{p,j}^T \Delta {\bf u},
\end{array}
\end{equation}

\noindent with $\Delta {\bf u} = (|u_1-u_{k,1} |, |u_2 - u_{k,2} |, \hdots, |u_{n_u} - u_{k,n_u}|)$. Defining the polytope

$$
\mathcal{\tilde L}_j = \left\{ {\bf u} : g_{p,j} ({\bf u}_k) + \tilde {\boldsymbol \kappa}_{p,j}^T \Delta {\bf u} \leq 0 \right\},
$$

\noindent we see that $\mathcal {\tilde L}_j \subseteq \mathcal{L}_j$ by virtue of (\ref{eq:proofineq}).

Consider now the Cauchy-Schwarz inequality, which states that

\begin{equation}\label{eq:cauchy}
\tilde {\boldsymbol \kappa}_{p,j}^T \Delta {\bf u} \leq \| \tilde {\boldsymbol \kappa}_{p,j} \|_2 \| \Delta {\bf u} \|_2 = \| \tilde {\boldsymbol \kappa}_{p,j} \|_2 \| {\bf u} - {\bf u}_k \|_2,
\end{equation}

\noindent and define the next (ball) subset as

$$
\mathcal{\hat L}_j = \left\{ {\bf u} : g_{p,j} ({\bf u}_k) + \| \tilde {\boldsymbol \kappa}_{p,j} \|_2 \| {\bf u} - {\bf u}_k \|_2 \leq 0 \right\}.
$$

\noindent That $\mathcal{\hat L}_j \subseteq \mathcal{\tilde L}_j$ follows from (\ref{eq:cauchy}). It now remains to show that $\mathcal{B}_e \subseteq \mathcal{\hat L}_j$, as this implies $\mathcal{B}_e \subseteq \mathcal{L}_j$ and completes the proof. For every ${\bf u} \in \mathcal{B}_e$, we have, by definition, $\| {\bf u} - {\bf u}_k \|_2 \leq \delta_e$, and so it follows that any ${\bf u} \in \mathcal{B}_e$ will also be in $\mathcal{\hat L}_j$ if

$$
g_{p,j} ({\bf u}_k) + \delta_e \| \tilde {\boldsymbol \kappa}_{p,j} \|_2 \leq 0
$$

\noindent holds, which is precisely what is ensured by (\ref{eq:backoff2}). \qed

\end{pf}

Employing the above theorem and corollary, we may now state the alternate results of Sections \ref{sec:consat}-\ref{sec:filter}:

\;
\noindent {\bf Constraint Satisfaction During Optimization (Section \ref{sec:consat})}
\;

\begin{equation}\label{eq:feasconddirlocC}
\hspace{-9mm}\begin{array}{l}
\displaystyle g_{p,j} ({\bf u}_k) + \sum_{i \in I_{{\rm ccv},j}^{{\bf u}_k,{\bf u}_{k+1}}} \mathop {\max} \left[ \begin{array}{l} \displaystyle \nabla \underline g_{p,ji}^{{\bf u}_k} ( u_{k+1,i} - u_{k,i} ), \vspace{1mm} \\ \displaystyle \nabla \overline g_{p,ji}^{{\bf u}_k} ( u_{k+1,i} - u_{k,i} ) \end{array} \right] \vspace{2mm} \\
\hspace{12mm}\displaystyle + \sum_{i \not \in I_{{\rm ccv},j}^{{\bf u}_k,{\bf u}_{k+1}}} \mathop {\max} \left[ \begin{array}{l} \underline \kappa_{p,ji}^{{\bf u}_k, {\bf u}_{k+1}} ( u_{k+1,i} - u_{k,i} ), \vspace{1mm} \\ \overline \kappa_{p,ji}^{{\bf u}_k, {\bf u}_{k+1}} ( u_{k+1,i} - u_{k,i} ) \end{array} \right] \leq 0 \Rightarrow g_{p,j} ({\bf u}_{k+1}) \leq 0.
\end{array}
\end{equation}

\;
\noindent {\bf Constraint Satisfaction During Perturbation (Section \ref{sec:excite})}
\;

\begin{equation}\label{eq:backoff3}
g_{p,j} ({\bf u}_k) + \delta_e \| \tilde {\boldsymbol \kappa}_{p,j} \|_2 \leq 0 \Rightarrow g_{p,j} ({\bf u}) \leq 0,  \; \forall {\bf u} \in \mathcal{B}_e.
\end{equation}

\;
\noindent {\bf Bounding Values in the Presence of Uncertainty (Section \ref{sec:filter})}
\;

\begin{equation}\label{eq:lipgenrobLdirlocC}
\underline f_p ^{\bar k} := \mathop {\max} \limits_{{\tilde k} = 0,...,k} \left( \begin{array}{l} \displaystyle \underline f_p^{\tilde k} + \sum_{i \in I_{\rm cvx}^{{\bf u}_{\tilde k}, {\bf u}_{\bar k}}} \mathop {\min} \left[ \begin{array}{l} \displaystyle \nabla \underline f_{p,i}^{\tilde k} ( u_{{\bar k},i} - u_{\tilde k,i} ), \vspace{1mm} \\ \displaystyle \nabla \overline f_{p,i}^{\tilde k} ( u_{{\bar k},i} - u_{\tilde k,i} ) \end{array} \right] \vspace{2mm} \\ \hspace{5mm} \displaystyle + \sum_{i \not \in I_{\rm cvx}^{{\bf u}_{\tilde k}, {\bf u}_{\bar k}}} \mathop {\min} \left[ \begin{array}{l} \underline \kappa_{i}^{{\bf u}_{\tilde k}, {\bf u}_{{\bar k}}} ( u_{{\bar k},i} - u_{\tilde k,i} ), \vspace{1mm} \\ \overline \kappa_{i}^{{\bf u}_{\tilde k}, {\bf u}_{{\bar k}}} ( u_{{\bar k},i} - u_{\tilde k,i} ) \end{array} \right]   \end{array} \right) \mathop \leq \limits^P f_p ({\bf u}_{\bar k}),  
\end{equation}

\begin{equation}\label{eq:lipgenrobUdirlocC}
f_p ({\bf u}_{\bar k}) \mathop \leq \limits^P \mathop {\min} \limits_{{\tilde k} = 0,...,k} \left( \begin{array}{l} \displaystyle \overline f_p^{\tilde k} + \sum_{i \in I_{\rm ccv}^{{\bf u}_{\tilde k}, {\bf u}_{\bar k}}} \mathop {\max} \left[ \begin{array}{l} \displaystyle \nabla \underline f_{p,i}^{\tilde k} ( u_{{\bar k},i} - u_{\tilde k,i} ), \vspace{1mm} \\ \displaystyle \nabla \overline f_{p,i}^{\tilde k} ( u_{{\bar k},i} - u_{\tilde k,i} ) \end{array} \right] \vspace{2mm} \\ \hspace{5mm} \displaystyle + \sum_{i \not \in I_{\rm ccv}^{{\bf u}_{\tilde k}, {\bf u}_{\bar k}}} \mathop {\max} \left[ \begin{array}{l} \underline \kappa_{i}^{{\bf u}_{\tilde k}, {\bf u}_{{\bar k}}} ( u_{{\bar k},i} - u_{\tilde k,i} ), \vspace{1mm} \\ \overline \kappa_{i}^{{\bf u}_{\tilde k}, {\bf u}_{{\bar k}}} ( u_{{\bar k},i} - u_{\tilde k,i} ) \end{array} \right]   \end{array} \right) := \overline f_p^{\bar k}.  
\end{equation}

Finally, it is worthwhile to show how the simpler bound of (\ref{eq:lipgenLU}) may be obtained from the alternate (\ref{eq:lipgenL})-(\ref{eq:lipgenU}). In doing so, we also present a sufficient condition for which (\ref{eq:lipgenL})-(\ref{eq:lipgenU}) provide tighter bounds over $\mathcal{I}_{{\bf u}_a}^{{\bf u}_b}$.

\begin{corol}[Relation Between (\ref{eq:lipgenLU}) and (\ref{eq:lipgenL})-(\ref{eq:lipgenU})]
Let $\mathcal{I}_{{\bf u}_a}^{{\bf u}_b} \subseteq \mathcal{I}$ and let $f_p$ be $\mathcal{C}^1$ over an open set containing $\mathcal{I}_{{\bf u}_a}^{{\bf u}_b}$, with its derivatives bounded as in (\ref{eq:liplocal})-(\ref{eq:derbounds}). Defining the mixed vectors of greatest absolute Lipschitz constants and derivative bounds as

$$
\tilde \kappa_{\overline f, i} = \left\{ \begin{array}{ll} \displaystyle \mathop {\max} \left( \Big | \nabla \underline f_{p,i}^{{\bf u}_a} \Big |, \Big | \nabla \overline f_{p,i}^{{\bf u}_a} \Big | \right), & i \in I_{\rm ccv}^{{\bf u}_a,{\bf u}_b}, \vspace{1mm} \\ \displaystyle \mathop {\max} \left( \Big |\underline \kappa_{i}^{{\bf u}_a,{\bf u}_b} \Big |, \Big |\overline \kappa_{i}^{{\bf u}_a,{\bf u}_b} \Big | \right), & i \not \in I_{\rm ccv}^{{\bf u}_a,{\bf u}_b}, \end{array} \right .
$$

$$
\tilde \kappa_{\underline f, i} = \left\{ \begin{array}{ll} \displaystyle \mathop {\max} \left( \Big | \nabla \underline f_{p,i}^{{\bf u}_a} \Big |, \Big | \nabla \overline f_{p,i}^{{\bf u}_a} \Big | \right), & i \in I_{\rm cvx}^{{\bf u}_a,{\bf u}_b}, \vspace{1mm} \\ \displaystyle \mathop {\max} \left( \Big |\underline \kappa_{i}^{{\bf u}_a,{\bf u}_b} \Big |, \Big |\overline \kappa_{i}^{{\bf u}_a,{\bf u}_b} \Big | \right), & i \not \in I_{\rm cvx}^{{\bf u}_a,{\bf u}_b}, \end{array} \right .
$$

\noindent it follows that $\| \tilde {\boldsymbol \kappa}_{\overline f} \|_2 \leq \kappa$ implies that (\ref{eq:lipgenU}) is tighter than the upper bound of (\ref{eq:lipgenLU}). Likewise, $\| \tilde {\boldsymbol \kappa}_{\underline f} \|_2 \leq \kappa$ implies that (\ref{eq:lipgenL}) is tighter than the lower bound of (\ref{eq:lipgenLU}).

\end{corol}
\begin{pf}
For the upper bound, we borrow the analysis already carried out in the proof of Corollary \ref{corol:backoff} to state that the bound

$$
f_p ({\bf u}_b) \leq f_p ({\bf u}_a) + \| \tilde {\boldsymbol \kappa}_{\overline f} \|_2 \| {\bf u}_b - {\bf u}_a \|_2
$$

\noindent is a looser, more conservative bound than (\ref{eq:lipgenU}). Clearly, $\| \tilde {\boldsymbol \kappa}_{\overline f} \|_2 \leq \kappa$ then implies that the upper bound of (\ref{eq:lipgenLU}) is even more conservative, whence the proof of the first half of the result.

For the lower bound, we note that

$$
\mathop {\min} \left[ \begin{array}{l} \nabla \underline f_{p,i}^{{\bf u}_a} (u_{b,i} - u_{a,i}), \vspace{1mm} \\  \nabla \overline f_{p,i}^{{\bf u}_a} (u_{b,i} - u_{a,i}) \end{array} \right] \geq -\mathop {\max} \left( \Big | \nabla \underline f_{p,i}^{{\bf u}_a} \Big |, \Big | \nabla \overline f_{p,i}^{{\bf u}_a} \Big | \right) | u_{b,i} - u_{a,i} |,
$$

$$
\mathop {\min} \left[ \begin{array}{l} \underline \kappa_{i}^{{\bf u}_a,{\bf u}_b} (u_{b,i} - u_{a,i}), \vspace{1mm} \\  \overline \kappa_{i}^{{\bf u}_a,{\bf u}_b} (u_{b,i} - u_{a,i}) \end{array} \right] \geq -\mathop {\max} \left( \Big |\underline \kappa_{i}^{{\bf u}_a,{\bf u}_b} \Big |, \Big |\overline \kappa_{i}^{{\bf u}_a,{\bf u}_b} \Big | \right) | u_{b,i} - u_{a,i} |,
$$

\noindent which allows us to apply the same analysis as before to obtain

$$
f_p ({\bf u}_b) \geq f_p ({\bf u}_a) - \| \tilde {\boldsymbol \kappa}_{\underline f} \|_2 \| {\bf u}_b - {\bf u}_a \|_2
$$

\noindent as a looser, more conservative alternative to (\ref{eq:lipgenL}). By the same reasoning, $\| \tilde {\boldsymbol \kappa}_{\underline f} \|_2 \leq \kappa$ then implies that the lower bound of (\ref{eq:lipgenLU}) is even more conservative, thus implying that (\ref{eq:lipgenL}) is tighter. \qed

\end{pf}

Clearly, we can obtain a valid version of (\ref{eq:lipgenLU}) from (\ref{eq:lipgenL})-(\ref{eq:lipgenU}) by setting $\kappa := \mathop {\max} (\| \tilde {\boldsymbol \kappa}_{\underline f} \|_2, \| \tilde {\boldsymbol \kappa}_{\overline f} \|_2)$.

\section{The Setting and Estimation of Lipschitz Constants}
\label{sec:estim}

Setting the Lipschitz constants in implementation so that they satisfy (\ref{eq:lipgen}) or (\ref{eq:liplocal}) is usually not trivial and requires some care. As already mentioned, choosing constants that do not satisfy these inequalities potentially makes the aforementioned techniques invalid because of the possibility that the corresponding Lipschitz bounds no longer hold. At the same time, as discussed at the start of Section \ref{sec:refine}, satisfying the inequalities with too much conservatism may affect performance undesirably.

It is thus the goal of this section to offer some insight into how one may set and refine the estimates of the Lipschitz constants intelligently. We begin with the special and desirable case where these constants are easily known from the physical laws governing a given experimental system (Section \ref{sec:physics}). Unfortunately, as there exists no means to rigorously guarantee that the Lipschitz constants chosen for a \emph{general experimental function} be correct,  one must inevitably turn to heuristic methods in the general case, and some of these are discussed in Sections \ref{sec:modelest}-\ref{sec:lipconsist}. Despite lacking the desired rigor, the methods proposed are not doomed to failure for this reason, and can certainly be applied with success, as has been demonstrated with some empirical evidence from both simulated and experimental case studies. This last point is discussed in Section~\ref{sec:empirical}.

\subsection{Setting Lipschitz Constants by Exploiting Physical Laws}
\label{sec:physics}

Virtually all experimental relationships are subject to physical laws, and for some relationships such laws are so well known and documented that they may be assumed to hold. A simple example previously employed by the authors is the exchange of heat between a heating and heated element -- e.g., a jacket (heating element) surrounding a reactor (heated element). If the temperature of the heated element is an experimental function, and the temperature of the heating element is one of the decision variables under the user's control, then much can be said about the sensitivity relating the changes in the temperatures of the two. For very many systems, it would be expected that raising the temperature of the heating element will raise the temperature of the heated element, and that, symmetrically, lowering the temperature of one will lead to a lower temperature in the other. As such, the sensitivity would always be positive, and the lower Lipschitz constant, $\underline \kappa_i$, could safely be set to 0. Additionally, heat losses are likely to ensure that a change in the temperature of the heating element will lead to a \emph{smaller} change in the temperature in the heated one, thus allowing us to set $\overline \kappa_i$ as $1$.

While this particular example is somewhat trivial, there are others that are less so. For example, in the solid oxide fuel-cell system where one manipulates the current while needing to respect a lower-limit constraint on the cell voltage \citep{Marchetti2009}, it is known that the relationship between voltage and current is always inverse proportional, thus allowing one to set the corresponding $\overline \kappa_i$ as 0. And although the physical laws could not tell us the value of $\underline \kappa_i$ in this case, one could nevertheless rely on the multitude of available experimental data to come up with a good estimate. For the polymerization example of Section \ref{sec:poly}, one may know that extended heating (increased $u_1$), together with shortened cooling (decreased $u_2$) and thus increased temperature, will generally lead to a favorization of the termination reactions and shorter polymer chains, thus resulting in a final product with a lower average molecular weight. This knowledge may then be used to set the Lipschitz constants for the molecular-weight constraint function accordingly.

Identifying such relationships and their relevance to the setting of Lipschitz constants may be an important step to include in the experimental and theoretical work that normally goes into understanding a given system prior to optimization. Note, however, that the physical laws are harder to link to the lumped $\kappa$ in (\ref{eq:lipgen}), which is more abstract than the individual $\underline \kappa_i$, $\overline \kappa_i$ of (\ref{eq:liplocal}).

\subsection{Model-Based Estimation}
\label{sec:modelest}

Many of the relationships that appear in experimental optimization problems will have undergone a fair amount of theoretical investigation, resulting in first-principles parametric models,

$$
f_p ({\bf u}) \approx f_{\hat p} ({\bf u}, {\boldsymbol \theta}),
$$

\noindent being available for them, where the uncertain parameters $\boldsymbol{\theta}$ are assumed to belong to a bounded set, $\Theta$.

If the assumption of parametric uncertainty is not too egregious, and if the uncertainty set $\Theta$ is not too erroneously defined, then a reasonable choice of Lipschitz constants in (\ref{eq:liplocal}) may be obtained by minimizing and maximizing the derivatives of the model over $\mathcal{I}_{{\bf u}_a}^{{\bf u}_b}$ and $\Theta$:

\begin{equation}\label{eq:lipmin}
\underline \kappa_{i}^{{\bf u}_a, {\bf u}_b} := {\rm arg} \mathop {{\rm{minimize}}}\limits_{\small \begin{array}{c} {\bf u} \in \mathcal{I}_{{\bf u}_a}^{{\bf u}_b} \\ {\boldsymbol \theta} \in \Theta \end{array}} \frac{\partial f_{\hat p}}{\partial u_i} \Big |_{{\bf u},{\boldsymbol \theta}},
\end{equation}

\begin{equation}\label{eq:lipmax}
\overline \kappa_{i}^{{\bf u}_a, {\bf u}_b} := {\rm arg} \mathop {{\rm{maximize}}}\limits_{\small \begin{array}{c} {\bf u} \in \mathcal{I}_{{\bf u}_a}^{{\bf u}_b} \\ {\boldsymbol \theta} \in \Theta \end{array}} \frac{\partial f_{\hat p}}{\partial u_i} \Big |_{{\bf u},{\boldsymbol \theta}},
\end{equation}

\noindent with the analogous formulation

\begin{equation}\label{eq:lipmax2}
\kappa := {\rm arg} \mathop {{\rm{maximize}}}\limits_{\small \begin{array}{c} {\bf u}_a, {\bf u}_b \in \mathcal{I} \\ {\bf u}_a \neq {\bf u}_b \\ {\boldsymbol \theta} \in \Theta \end{array}} \frac{f_{\hat p}({\bf u}_b, {\boldsymbol \theta}) - f_{\hat p}({\bf u}_a, {\boldsymbol \theta})}{\| {\bf u}_b - {\bf u}_a \|_2}
\end{equation}

\noindent for the lumped $\kappa$. While solving these optimization problems may not be computationally easy, it is conceivable that they would be solved off-line prior to any $\Gamma$ being applied, thus making the associated computational burdens of lesser concern. If an additional layer of safety is needed, one can heuristically lower and increase the estimates as required.

For the very general case where no model is available, one could still apply the above approach but in a data-driven fashion. One could, for example, construct a linear or quadratic model from several obtained measurements, choose ${\boldsymbol \theta}$ as the coefficients of the model, let their confidence intervals define $\Theta$ \citep[\S 10.5]{Montgomery2012}, and then apply either (\ref{eq:lipmin})-(\ref{eq:lipmax}) or (\ref{eq:lipmax2}).

\subsection{Refinement via Data-Driven Consistency Checks}
\label{sec:lipconsist}

The previous two sections have addressed the problem of setting Lipschitz constants when none are initially provided. Now, let us consider a technique for refining these initial estimates so that they are \emph{consistent with the obtained data}.

Using (\ref{eq:lipgen}) as an example, suppose now that we have obtained an estimate of a Lipschitz constant, denoted by $\hat \kappa$, and as such have the hypothetically valid bound

\begin{equation}\label{eq:lipboundest}
f_{p} ({\bf u}_b) \leq f_{p} ({\bf u}_a) + \hat \kappa \| {\bf u}_b - {\bf u}_a \|_2, \; \forall {\bf u}_a, {\bf u}_b \in \mathcal{I}.
\end{equation}

\noindent Although there exists no way to prove the validity of (\ref{eq:lipboundest}) for an experimental function $f_{p}$ in the general case, we can and should confirm that this bound is at least satisfied for those experimental measurements that have been collected.

Letting $k_1$ and $k_2$ denote two (different) indices, a basic consistency-check algorithm would consist in the following two steps being carried out for every combination $(k_1,k_2)$ of  $k_1, k_2 \in [0,...,k]$:

\begin{enumerate}[(i)]
\item Check if the inequality

\begin{equation}\label{eq:lipboundest2}
f_{p} ({\bf u}_{k_2}) \leq f_{p} ({\bf u}_{k_1}) + \hat \kappa \| {\bf u}_{k_2} - {\bf u}_{k_1} \|_2
\end{equation}

\noindent is satisfied.

\item If (\ref{eq:lipboundest2}) is satisfied, proceed to the next $(k_1, k_2)$ combination unless no combinations are left, in which case terminate. Otherwise, increase the value of $\hat \kappa$ by a preset, strictly positive quantity and return to (i).

\end{enumerate}

\noindent It is easy to show that such an algorithm will terminate with a $\hat \kappa$ value that is consistent with the obtained data. For the more complex case of (\ref{eq:liplocal}), one could apply the same concept by lowering $\underline \kappa_i$ and raising $\overline \kappa_i$ until (\ref{eq:lipgenL})-(\ref{eq:lipgenU}) hold. While there is no apparent ``best way'' to increase the $\underline \kappa_{i}$, $\overline \kappa_i$ values when inconsistency occurs, one could employ heuristics, such as doubling the values. More elegant schemes may also be proposed \citep[\S 3.5.3]{Bunin:SCFOImp}.

\subsection{Empirical Evidence for Lipschitz-Constant Estimation}
\label{sec:empirical}

Clearly, since Lipschitz constants are often unknown, it is difficult to test, in any sort of closed-loop manner, whether different estimation methods are actually successful at identifying proper values for these constants in practice. Instead, one may apply the methods during the solution procedure of (\ref{eq:mainprob}) and see if the problem is solved in a satisfactory manner. When the problem comes with experimental constraints, the strongest criterion is to verify that the constraints were met during optimization -- if they were not, then this is strong evidence for poor estimates, as it suggests that the constraint-satisfaction condition (\ref{eq:congamma}) was not valid for the constraints violated. The speed of decrease in the cost function value is another telling criterion -- if progress is very slow, then this may be a sign of the estimates being too conservative. In other cases, the quality of the estimates may be difficult to evaluate, although different ``strange'' behaviors during optimization may be an indicator that the constants are not being estimated properly. To date, much of the available empirical evidence for the effectiveness of Lipschitz-constant estimation methods comes from problems that have been tackled by the SCFO solver \citep{SCFOug}, as this solver incorporates almost all of the techniques discussed in this paper and relies heavily on the Lipschitz constants to operate. 

Of the experimental results available, perhaps the most telling are those obtained in the experimental optimization of a laboratory solid-oxide fuel cell stack \citep{Francois2014}, where no theoretical model was used and the Lipschitz constants were initialized using a data-driven model and the methods of Section \ref{sec:modelest}, before being refined by the consistency checks of Section \ref{sec:lipconsist} as more data became available. From the limited results obtained, one does see that the estimation is sufficiently good from the safety perspective -- the constants are conservative enough so as to keep the system from violating its two experimental constraints, and a close-to-optimal system efficiency is achieved without requiring too many iterations, which suggests that the constants are not \emph{too} conservative.

Other experimental results have been obtained by applying an earlier version of the SCFO to the problem of run-to-run controller tuning, where the tuning parameters of the controllers were treated as the decision variables, with different trajectory tracking metrics defined to judge controller performance for a given setpoint profile \citep{Bunin2013p}. The solver managed to successfully autotune a model-predictive controller for a water-tank system and a fixed-order controller for a mechanical torsional plant -- however, these results are less validating since no experimental constraints were present and so it was impossible to evaluate the quality of the Lipschitz-constant estimates with regard to their ability to enforce constraints. At the same time, the constants for the experimental cost function \emph{were} used to help reduce the effects of variance in the cost-function values (Section \ref{sec:filter}). While this is likely to have aided in the optimization, there is no way to say how much, as testing the quality of the Lipschitz-constant estimators was not the goal of this largely proof-of-concept study. Once again, no \emph{a priori} models were used -- the solver constructed a data-driven model from $n_u+1$ initial experiments, and then refined the initial estimates via the consistency checks.

In simulated tests, a good quantity of results is now available for a number of case-study test problems, including six with experimental constraints, courtesy of the ExpOpt database \citep{ExpOpt}. Examining the results of the SCFO solver for these problems shows that the current estimation methods employed by the solver are certainly not perfect, and one can indeed have Lipschitz-constant estimates that are not valid and that allow for constraint violations, although this varies significantly from problem to problem. For certain problems, the estimates are consistently good and violations are minimal, with the violations being of very small magnitude even when they are present. For other problems, violations may be larger and more frequent, but it is nevertheless rare to see consistently large constraint violations, suggesting that consistency checks are eventually able to improve on poor initial estimates.

With regard to results obtained for methods other than the SCFO, an application of a data-driven approach of Section \ref{sec:modelest} has recently been reported by \cite{Bunin2016}. Here, a modified version of G. E. P. Box's evolutionary-operation method \citep{Box:69} was proposed, with the Lipschitz bound used to ensure that the experiments carried out by the algorithm would never violate the problem constraints. After each set of $2n_u$ axially distributed experiments, the Lipschitz constants were estimated locally as the derivatives of a linear model for the $2n_u+1$ axial-plus-center data points, with additional conservatism added to the estimates to account for noise and nonlinearity effects. The end result was very satisfactory -- in testing the algorithm for three problems from the ExpOpt database, the estimated Lipschitz constants were always sufficiently conservative and the experimental constraints were not violated once.

\section{Concluding Remarks}
\label{sec:conclude}

It is the authors' hope that the present document has helped convince the reader that Lipschitz constants -- despite being an implicit, age-old mathematical concept -- may be made explicit and thereby bring multiple benefits to the experimental optimization domain. Most importantly, the use of these constants appears to offer an extremely simple way to guarantee constraint satisfaction in problems with experimental constraints, which is something that may be exploited during both the optimization and information-gathering (perturbation) phases of the solution process. Because constraint satisfaction may be immensely important in experimental settings, this alone presents a very strong argument for the implementation of Lipschitz constants and bounds. Additionally, it has been shown that these constants may play a useful role in reducing the effect of measurement or estimation uncertainty, which, depending on the optimization algorithm, may lead to significant performance improvements. There are no major implementation difficulties involved with the techniques described, and so it is the authors' recommendation that these techniques be given consideration when approaching experimental optimization problems.

This paper has also addressed the two major concerns regarding the use of Lipschitz constants in practice -- the concern that the set constants may be too conservative and the concern that they may be improperly set or estimated. Both are valid and represent potential issues. However, as has been shown in Section \ref{sec:refine}, there do exist several ways to reduce conservatism -- namely, one may tighten the bounds by exploiting direction, locality, and the potential convexity/concavity properties of the function. With regard to estimation, Section \ref{sec:estim} has outlined a number of methods to set and refine the Lipschitz-constant estimates, and empirical evidence has shown that such methods can, in fact, be sufficient for many problems.

\section*{Acknowledgements}
 
The authors would like to thank Professor Dominique Bonvin of the Laboratoire d'Automatique (\'Ecole Polytechnique F\'ed\'erale de Lausanne) for his insights and input.

%\bibliography{Lip}             % bib file to produce the bibliography
                                                    % with bibtex (preferred)

\end{document}